\documentclass{amsart}
\usepackage{amsmath,amsfonts,amssymb,enumerate,latexsym,graphicx,pstricks}
\input{epsf.tex}

\numberwithin{equation}{section}

\theoremstyle{plain}
\newtheorem{theorem}{Theorem}[section]
\newtheorem{corollary}[theorem]{Corollary}
\newtheorem{proposition}[theorem]{Proposition}

\newtheorem{lemma}{Lemma}[section]
\newtheorem*{hypothesis}{Inductive Hypothesis(n)}

\newtheorem*{mainlemma}{Main Lemma 5.6}

\theoremstyle{definition}
\newtheorem{definition}{Definition}[section]

\newtheorem{numberedremark}{Remark}[section]

\theoremstyle{remark}

\begin{document}
\title{Convex Dynamics and Applications}
%\today\ draft

\author[R. L. Adler]{R. L. Adler}
\address{R.L. Adler\\IBM, TJ Watson Research Center\\
Yorktown Heights, NY 10598-0218\\
USA} \email{rla@us.ibm.com}

\author{B. Kitchens}
\address{B. Kithchens\\Mathematics Department\\
State University of New York at Stony Brook\\
Stony Brook\\
NY 11794\\
USA} \email{kitchens@math.sunysb.edu}

\author[R. Martens]{M. Martens}
\address{M. Martens\\University of Groningen\\
Department of Mathematics\\
P.O.Box 800\\
9700 AV Groningen\\
The Netherlands} \email{marco@math.rug.nl}

\author[C. Pugh]{C. Pugh}
\address{C. Pugh\\Mathematics Department\\
University of California \\
Berkeley California, 94720, \\
USA} \email{pugh@math.berkeley.edu}

\author[M. Shub]{M. Shub}
\address{M. Shub\\
IBM, TJ Watson Research Center\\
Yorktown Heights, NY 10598-0218\\
USA}
\curraddr{\\Department of Mathematics\\
University of Toronto\\
100 St. George Street\\
Toronto, Ontario M5S 3G3 Canada}
 \email{shub@math.toronto.edu}
\thanks{Shub's work was partly funded by NSF Grant \#DMS-9988809}
\thanks{Martens' and Tresser's work was partly funded by NSF Grant\#DMS-0073069}
\author[C. Tresser]{C. Tresser}
\address{C. Tresser\\IBM, TJ Watson Research Center\\
Yorktown Heights, NY 10598-0218\\
USA} \email{tresser@us.ibm.com}

\keywords{dynamical system, orbit bounding, convexity, greedy
algorithm, error diffusion}
\date{\today}
\begin{abstract}
This paper proves a theorem about bounding orbits of a time
dependent dynamical system. The maps that are involved are
examples in convex dynamics, by which we mean the dynamics of
piecewise isometries where the pieces are convex. The theorem came
to the attention of the authors in connection with the problem of
digital halftoning. \textit{Digital halftoning} is a family of
printing technologies for getting full color images from only a
few different colors deposited at dots all of the same size. The
simplest version consist in obtaining grey scale images from only
black and white dots. A corollary of the theorem is that for
\textit{error diffusion}, one of the methods of digital
halftoning, averages of colors of the printed dots converge to
averages of the colors taken from the same dots of the actual
images. Digital printing is a special case of a much wider class
of scheduling problems to which the theorem applies. Convex
dynamics has roots in classical areas of mathematics such as
symbolic dynamics, Diophantine approximation, and the theory of
uniform distributions.
\end{abstract}

\maketitle
\thispagestyle{empty} \def\IMSmarkvadjust{0 pt}
\def\IMSmarkhadjust{0 pt}
\def\IMSmarkhpadding{0 pt}
\def\IMSpubltext{Published in modified form:}
\def\SBIMSMark#1#2#3{
 \font\SBF=cmss10 at 10 true pt
 \font\SBI=cmssi10 at 10 true pt
 \setbox0=\hbox{\SBF \hbox to \IMSmarkhpadding{\relax}
                Stony Brook IMS Preprint \##1}
 \setbox2=\hbox to \wd0{\hfil \SBI #2}
 \setbox4=\hbox to \wd0{\hfil \SBI #3}
 \setbox6=\hbox to \wd0{\hss
             \vbox{\hsize=\wd0 \parskip=0pt \baselineskip=10 true pt
                   \copy0 \break%
                   \copy2 \break% 
                   \copy4 \break}}
 \dimen0=\ht6   \advance\dimen0 by \vsize \advance\dimen0 by 8 true pt
                \advance\dimen0 by -\pagetotal
	        \advance\dimen0 by \IMSmarkvadjust
 \dimen2=\hsize \advance\dimen2 by .25 true in
	        \advance\dimen2 by \IMSmarkhadjust

%
%   Check for publication info
%
%  \newread\jref
  \openin2=publishd.tex
  \ifeof2\setbox0=\hbox to 0pt{}
  \else 
     \setbox0=\hbox to 3.1 true in{
                \vbox to \ht6{\hsize=3 true in \parskip=0pt  \noindent  
                {\SBI \IMSpubltext}\hfil\break
                \input publishd.tex 
                \vfill}}
  \fi
  \closein2
  \ht0=0pt \dp0=0pt
 \ht6=0pt \dp6=0pt
 \setbox8=\vbox to \dimen0{\vfill \hbox to \dimen2{\copy0 \hss \copy6}}
 \ht8=0pt \dp8=0pt \wd8=0pt
 \copy8
 \message{*** Stony Brook IMS Preprint #1, #2. #3 ***}
}

\def\IMSmarkvadjust{-30pt}
\SBIMSMark{2004/02}{February 2004}{}

\bibliographystyle{ieeetr}

% Section 1
% Section 1
\section{Introduction}
% Section 1
% Section 1

\begin{theorem}\label{T:main}
Let $P$ be a polytope in $\mathbb{R}^N$, and for each $\gamma$ in
$P$ let $\phi_\gamma$ denote the map from $\mathbb{R}^N$ to
$\mathbb{R}^N$ defined by

\begin{equation}\label{E:map}
 \phi_{\gamma}(x)=x+(\gamma-v(x)),
\end{equation}
where $v(x)$ is the closest vertex of $P$ to $x$ in Euclidean norm
with some tie-breaking rule. Then for any compact set K there
exists a bounded convex domain $Q$ containing $K$ and invariant in
the sense that for all $\gamma$ in $P,\ \phi_\gamma Q \subset Q \,
.$

\end{theorem}

 We encountered Theorem \ref{T:main} in connection with a problem in digital
 color printing \cite{AKMTW2002}. It turns out to have wide applicability (see,
\emph{e.g.,} \cite{niederreiter:applications71},
\cite{GAGA:VBoundsOnSched:2003}) and considerable mathematical
depth.

\medskip
The general problem that this theorem addresses concerns
approximating a bounded one-sided sequence $\gamma(k)$ of
arbitrary values in $\mathbb{R}^N$ by a sequence of elements
$V(k)$ chosen from some finite set. In a wide range of
applications, such as digital printing or scheduling, a reasonable
sense of well-approximation is that there be a uniform bound
(small if possible) on some norm of the \textit{cumulative error
vectors}
\begin{equation}\label{E:error}
\varepsilon(n)=\sum_{k=0}^n(\gamma(k) - V(k)),
\end{equation}
in which case the average error will go to zero. In applications
it is assumed that the $V(k)$'s are taken as the vertices of a
polytope $P$ and the $\gamma(k)$'s belong to this polytope. In the
case when the norm $||\cdot||$ is the sup norm and $P$ the
standard simplex, the problem of finding the best universal bounds
on $\sup _n||\varepsilon(n)||$ is often referred to as the
\textit{Chairman Assignment Problem}. This problem was posed by
Niederreiter in \cite{Nied} and considered in
\cite{MeijNie,M73,T73,T80}; we say \textit{Generalized Chairman
Assignment Problem} when the norm is not necessarily the $\sup$
norm. Using the methods of \cite{GAGA:VBoundsOnSched:2003} to find
sequences that solve the Generalized Chairman Assignment Problem
turns out not to be effective as it requires looking ahead with
unbounded waiting time, at least for the $\sup$ norm. Waiting for
future input data before making output decisions is not practical
in many cases and sometimes even not acceptable as the scheduling
algorithms may need, not only to be \emph{causal} (\emph{i.e.,}
depend only on the past and the present and not also the future)
but also to be \emph{on the fly} (\emph{i.e.,} provide outputs
essentially as soon as input data are known). Thus, if one is not
looking for the best bound but only for a uniform one, or if one
needs to perform continual optimization, then a surprisingly good
approach is to write $\varepsilon(k)$ recursively as
\begin{equation}\label{recursion}
\varepsilon(k+1) =\varepsilon(k)+
\gamma(k+1)-v(\varepsilon(k)+\gamma(k+1))\,,
\end{equation}
and use the \textit{greedy algorithm} to specify
$v(\varepsilon(k)+\gamma(k+1))$ which means here that
$v(\varepsilon(k)+\gamma(k+1))$ is the vertex which minimizes the
norm of $\varepsilon(k+1)$ with some tie breaking rule:
\textit{i.e.,} the vertex closest to $\varepsilon(k)+\gamma(k+1)$
with some tie breaking rule. We shall be interested in the
Euclidean norm; and in what follows, when we say \emph{``the"
greedy algorithm}, we mean the greedy algorithm using the
Euclidean norm with any tie breaking rule for points equidistant
from two or more vertices. Later it will be convenient to express
the greedy algorithm in terms of the notion of Vorono\"{\i}
regions. A \emph{Vorono\"{\i} region} of a vertex is the closure
of the set of points that are closer to the vertex than to any
other vertex (see Definition \ref{D:Vor}).  Alternatively
expressed, the greedy algorithm chooses the vertex in the
Vorono\"{\i} region containing $\varepsilon(k)+\gamma(k+1)$ with
some decision rule when $\varepsilon(k)+\gamma(k+1)$ lies on a
boundary.

\medskip
%Recursion \eqref{recursion} defines a non-autonomous dynamical system of the kind encountered in control theory.
Recursion \eqref{recursion} defines a non-autonomous dynamical
system acting on the space of errors with a time-dependant
parameter that belongs to the polytope. We are going to focus our
attention on an associated dynamics in the ambient space of the
polytope. We make the following changes of variables:
\begin{equation}\label{E:changeofvariables}
x(k)=\gamma(k)+\varepsilon(k-1)\,,
\end{equation}
and
\begin{equation}\label{E:greedyvertex}
V(k)=v(x(k))=v(\gamma(k)+\varepsilon(k-1))\,.
\end{equation}
Then adding $\gamma(k+2)$ to both sides of \eqref{recursion} and
reducing indices, we get
\begin{equation}\label{recursion2}
x(k+1)=x(k)+\gamma(k+1)-v(x(k))\,.
\end{equation}
Recall that $v(x(k))$ is chosen as the closest vertex to $x(k)$
with some tie breaking rule. The orbit $x(k)$ can be expressed in
terms of the mapping $\phi_\gamma$ in \eqref{E:map} by
\begin{equation}\label{E:orbitmap}
x(k)=x(k-1)+(\gamma(k)-v(x(k-1)))=\phi_{\gamma(k)}(x(k-1))\,,
\end{equation}
and bounding it in Euclidean norm is answered by Theorem
\ref{T:main}. Bounding $\varepsilon(k)$ is equivalent to bounding
$x(k),$ and an immediate corollary of Theorem \ref{T:main} is the
following convergence of averages:
\begin{equation}\label{E:averages}
\lim_{n \to \infty}||\frac{1}{n}\sum_0^{n-1}
\gamma(k)-\frac{1}{n}\sum_0^{n-1} v(x(k))|| = 0\,.
\end{equation}

\subsection{The greedy algorithm in digital printing.}
To see the relevance of the above problem of approximating the
$\gamma(k)$'s by the $V(k)$'s to color printing requires some
background in color theory. The color printing alluded to concerns
\textit{digital halftoning} which is the printing technique for
imitating full color natural images (as opposed to colors placed
at random on a page) by a checkerboard of dots limited to a few
available colors. The simplest halftoning problem is that of
obtaining gray-scale images with only black and white dots.

\medskip
Color theory itself is a fascinating subject which we can only
briefly touch upon. It has long been known that our perception of
color has a vector space model. Electromagnetic spectra are so
rich that to describe them requires an infinite dimensional
function space. Nevertheless, experiments indicate that the human
visual system reduces the visible part of this space to a three
dimensional convex cone. In 1931 the \textit{Commission
Internationale de L'Eclairage} (CIE) devised a certain three
dimensional coordinate system called the \textit{tristimulus
space} which quantifies this cone. A complete account of how this
coordinate system is contrived from color matching experiments is
to be found in \cite{SW82}. Our perception of light color is
governed by a set of rules known as Grassman's laws, the
consequence of which is that two different linear combinations of
color vectors appear as the same color perhaps with different
radiances, if their sums are collinear. Consider a section of
paper on which a color image is to be printed using a digital
printer, typically a laser or ink jet printer which respectively
uses toner or ink to render colors. We shall side step the
nonmathematical practical difficulties involved in halftoning and
assume an idealized situation: namely, the paper upon which an
image is to be printed is partitioned into a checkerboard of tiny
squares called \textit{pixels} after a term used in color TV. We
shall label pixel locations as we would the entries of a matrix.
On each pixel is deposited a uniform color restricted to a very
limited number of available choices which we shall call
\textit{pixel colors}. The problem is how should one choose pixel
colors to render an acceptable full color image.

Over a small area that still consists of a multitude of pixels the
eye averages the received light. This means that we do not
perceive colors on a page by taking positive linear combinations
of color vectors, but rather convex combinations. The pixel colors
commonly available to a printer are the standard inks or toner
colors: cyan($C=(21,27,72)$), magenta($M=(33,18,22)$),
yellow($Y=(65,76,14)$). In addition red($R=(30,18,7)$),
green($G=(11,22,13)$), and blue($B=(9,7,20)$) can be gotten by
some kind of pairwise mixing of the above colors. Also available
are the ink or toner for black($K=(5,6,6)$), and
white($W=(84,87,105)$) which is taken for the color of paper in
the absence of any ink or toner. In parentheses are typical
tristimulus space coordinates for color vectors. The finite set
$\{ C,\ M,\ Y,\ R,\ G,\ B,\ W,\ B\}$ of ideally available printer
color vectors are vertices of a convex polytope, in fact
topologically a cube in three space. For a general reference about
how color theory affects digital color printing see \cite{K}.

\medskip
Since color usually varies with low spacial frequency over many
portions of natural images, it seems natural to simplify the
halftoning problem by asking how should one choose pixel colors to
represent a single uniform non-pixel color. Furthermore, since
averaging plays such an essential role in color perception, it
seems appropriate to recall that the uniform distribution theory
developed by T. van Aardenne-Ehrenfest, H. Davenport, K.F. Roth
and others. Their results concerning the discrepancy between
average values and expected ones can be applied to the problem of
approximating arbitrary uniform greys over two dimensional regions
by combinations of black and white squares on a lattice. From a
theorem of Roth (see \cite{BC}, p. 3) one can deduce that it is
not possible to bound the cumulative error simultaneously for all
grey levels and all rectangles. This is in contrast with printing
in one dimension--\textit{i.e.} printing on a line instead of a
page-- where it is well known, and as Equation \eqref{E:errorrun}
below shows, that the cumulative error can be very well bounded on
all segments, not only for grey levels but also for color. The
discussion so far seems to suggest that, while there is no
inherent difficulty in performing digital printing on a line,
there may be one a page, where it is of practical importance.
However, the issue is not one of bounding cumulative errors, but
rather of limiting their growth rate. Furthermore, a two
dimensional region $R$ of pixels over which the eye averages is
special.

\medskip
For the sake of simplification, let $R$ be a $n\times n$ square of
pixels. The local cumulative error $E(R)$ made over the region $R$
is given by
\begin{equation}\label{E:2derror}
E(R)=\sum_{(i,j) \in R} (\gamma(i,j)-V(i,j)).
\end{equation}
In the terminology of printing, the $\gamma(i,j)$'s are called
\textit{inputs} and the $V(i,j)$'s are called \textit{outputs}.
The object of halftoning is to chose outputs so that the average
cumulative error $|| E(R) || / n^2 $ is small. With respect to a
polytope $P$ in color space whose vertices are the inputs, as we
have seen, the greedy algorithm of Equation \eqref{recursion} will
achieve this. This is one of the methods of halftoning which is
called \textit{error diffusion} \cite{FS,JJN,S,K,Uli1987}. We call
this particular version \textit{simple error diffusion}
\cite{AKMTW2002}. Since this involves a sequential process, pixels
must linearly ordered, the usual order being the lexicographic
one: namely, for an $M\times N$ array the $(i,j)^{th}$ pixel gets
index $k=i\cdot N+j$ where the indexing starts at 0. If $B$ is a
bound on $\varepsilon(k),$ then $2B$ is a bound on the size of
cumulative error
\begin{equation}\label{E:errorrun}
\varepsilon(L+n-1)-\varepsilon(k)=\sum^{L+n-1}_{j=L}(\gamma(j)-V(j))
\end{equation}
due a run of $n$ pixels starting at pixel $L.$. Since there are
$n$ such runs in $R,$ we have that $||E(R)|| \le 2Bn$ which leads
to small average size of cumulative error for large enough $n$.

For two dimensional digital printing simple error diffusion is
deficient because according to Equation \eqref{recursion}, even
forgetting about the difficulty associated with changing lines,
only pixel $k-1$ has a direct influence on pixel $k.$ Other nearby
pixels have no direct influence. Even using more exotic orderings
cannot eliminate this flaw. However, one can simultaneously
improve the influence nearby pixels and deal with edge effects by
using a \textit{general error diffusion} which involves the greedy
algorithm with respect to a weighting of past errors as follows.
Let $w^k_i, 1\le i\le m,$ be a system of weighting factors chosen
in most cases for each $k$ to be a probability vector. Indeed we
shall always assume that the weighting factors are chosen this
way: \textit{i.e.,}
\[
w^k_i\ge 0, \text{ and } \sum_{i=1}^m w^k_i =1 \,.
\]
Define
\begin{equation}\label{E:pasterrors}
\varepsilon(k+1)=\left[\sum_{i=1}^m w_i^k
\varepsilon(k+1-i)\right] +\gamma(k+1)-V(k+1)
\end{equation}
where $V(k+1)$ minimizes $||\varepsilon(k+1)||$ with some
tie-breaking rule.
\begin{table}
 \begin{center}
  \begin{tabular}{| l | l | l | l | l | l |}
   \hline
    $w_m^k=p_{12}$  &$w_{m-1}^k=p_{11}$   &$w_{m-2}^k=p_{10}$
                          &$w_{m-3}^k=p_{9}$   &$w_{m-4}^k=p_{8}$  \\ \hline
    $w_{\ell-4}^{k}=p_{7}$  &$w_{\ell-3}^{k}=p_{6}$   &$w_{\ell-2}^{k}=p_{5}$
                          &$w_{\ell-1}^{k}=p_{4}$   &$w_{\ell}^{k}=p_{3}$  \\ \hline
    $w_{2}^k=p_{2}$ &$w_{1}^k=p_{1}$ &k+1\textsuperscript{st} pixel
                          &\multicolumn{2}{l |}{\dots \  $\ell=k+3-N$}   \\
   \hline
\end{tabular}
\smallskip
 \caption{Table of nonzero weights}\label{Ta:1}
\end{center}
\end{table}
In the lexicographic ordering of pixels, even though those to the
right and below have no influence on the current pixel, this
method of weighting the past orbit seems to be a much better
approximation to the eye's method of averaging than that of giving
all the weight to the pixel immediately preceding the current one.
A typical scheme for weighting factors for pixels not near borders
is shown in Table \eqref{Ta:1} where $p_i$ is a twelve-dimensional
probability vector. Since this vector is allowed to depend on $k,$
adjustments can be made to deal with pixels near borders.

The first such weighting scheme was introduced in \cite{FS} by R.
Floyd and L. Steinberg in the seminal paper where they first
described diffusion as a method for digital printing in 1975; this
original scheme only involved $p_1,\ p_5,$ and $p_6,$ all other
weighting factors being 0. Two more elaborate schemes can be found
in \cite{JJN} and \cite{Stuck} using the full Table \ref{Ta:1} and
a lot more have been devised and utilized since. In these two
schemes, as is often the case now for square pixels, the highest
weights are given to $p_1=p_5$, the next higher to $p_4=p_6$ and
so on , with value decreasing with the distance to the $k^{\rm
th}$ pixel being treated (for specific values and other
information concerning digital halftoning see the above references
as well as \cite{K} and \cite{Uli1987}).

\medskip

For general error diffusion the relation of $\varepsilon(k)$ to
cumulative error is more complicated than Equation
\eqref{E:errorrun}. Nevertheless we shall show that, as in the
case of simple error diffusion, for general error diffusion
$||E(R)||/n^2\rightarrow 0.$

\medskip
First we have the following corollary to Theorem \ref{T:main}.
\begin{corollary}\label{C:main} For general error diffusion
$||\varepsilon(k)||<B$ for some positive number $B.$
\end{corollary}

\begin{proof}
 Setting
\begin{equation}\label{E:changeofvariables2}
x(k+1)=\left[\sum_{i=1}^m w_i^k \varepsilon(k+1-i)\right]
+\gamma(k+1)
\end{equation}
(in the literature this is called the \textit{modified input})we
have
\[
\varepsilon(k)=x(k)-v(x(k)),
\]
where v(x(k))=V(k) is the closest vertex to x(k) with some
tie-breaking rule. Then
\begin{align}\label{E:weightedpastorbit}
x(k+1)&=\sum_{i=1}^m w_i^k
\left[x(k+1-i)-v(x(k))\right]+\gamma(k+1)\\
      &=\sum_{i=1}^m w_i^k\phi_{\gamma(k+1)}(x(k+1-i))\,.
\end{align}
The fact that $x(k)$ of \eqref{E:weightedpastorbit} is bounded,
and equivalently that $\epsilon(k)$ is bounded, is a consequence
of the convexity of Q in Theorem \eqref{T:main}.$\Box$
\end{proof}

\medskip
It is folklore knowledge in digital printing that for an $n\times
n$ square $R$, the average error $||E(R)||/n^2$ goes to zero with
$n$ (see for example \cite{Broj,Ko}). For the sake of completeness
we indicate a proof: more precisely that $||E(R)||=O(n)$. Suppose
that the pixel locations of the $n\times n$ square $R$ contained
in an $M\times N$ array, $N>>n,$ are $L+kK+j,\ 0\le j,k\le n-1.$
Then from Equation \eqref{E:pasterrors} we have
\begin{align}\label{E:2derror2}
   E(R)&=\sum_{k=0}^{(n-1)}\ \sum_{j=L+kN}^{L+kN+n-1}
        (\gamma(j)-V(j))\\
   &=\sum_{i=1}^m w_i\ \sum_{k=0}^{n-1}\ \sum_{j=L+kN}^{L+kN+n-1}\
       (\varepsilon(j) -\varepsilon(j-i)) \notag
\end{align}
The only nonzero weights are $\{ w_1, w_2, w_{pN+q}:\ 1 \le p \le
2,\ -2\le q \le 2 \}.$

Changing the limits on the index $j,$ and interchanging order of
summation \eqref{E:2derror2} we get
\begin{equation}\label{2dcumulative}
 E(R)= \sum_{i=1}^m w_i\  \sum_{j=0}^{n-1}\ \left[ \sum_{k=0}^{n-1}\
         \varepsilon(j+L+kN) - \sum_{k=0}^{n-1}\
         \epsilon(j+L+kN-i))\right]. \notag
\end{equation}
Let
\begin{equation}\label{E:E_i}
E_i\equiv   \sum_{j=0}^{n-1}\ \left[ \sum_{k=0}^{n-1}\
         \varepsilon(j+L+kN) - \sum_{k=0}^{n-1}\ \epsilon(j+L+kN-i)\right]
\end{equation}
Then for $0\le p\le 2, \ -2\le q \le 2,$
\begin{align}
E_{pN+q}&=  \sum_{j=0}^{n-1}\ \left[ \sum_{k=0}^{n-1}\
         \varepsilon(j+L+kN) - \sum_{k=0}^{n-1}\ \epsilon(j+L+(k-p)N-q)\right]\\
         &= \sum_{j=0}^{n-1}\ \sum_{k=0}^{n-1}\
         \varepsilon(j+L+kN) - \sum_{j=-q}^{n-q-1}\ \sum_{k=-p}^{n-p-1}\
         \epsilon(j+L+kN).\notag
\end{align}
We leave to the reader the subsequent calculations, the result of
which is
\[
 ||E(R)||=O(n).
 \]

\subsection{A pursuit problem.}
The greedy algorithm also leads to a solution to the following
pursuit problem.  Suppose a predator is chasing a prey in a
polytope, and the following conditions are imposed. The positions
of the predator and prey at time $t_n=t_{n-1}+\frac{1}{n}$ are
denoted by $p(n)$ and $q(n)$ respectively. We take $n$ to be
greater than zero and assume for that $p(0)$ and $q(0)$ are at two
distinct corners, but any $p(0)\neq q(0)$ would do. Between times
$t_n$ and time $t_{n+1}$ the prey moves in the direction of any
point $\gamma(n+1)$ in the polytope with speed $\frac{n}{n+1}$
while the predator moves at unit speed but restricted towards some
vertex, say $V(n+1),$ the choice of which being dictated by a
pursuit strategy. Notice that the prey has the advantage over the
predator of freer movement but the disadvantage of slower speed
(although less and less so). The rules of movement are expressed
by the following formulas:
\begin{align}
&p(n+1)=p(n)+ \frac{1}{n}(V(n+1)-p(n))\,,\label{E:target} \\
&q(n+1)=q(n)+\frac{1}{n+1}(\gamma(n+1)-q(n))\,.\label{E:pursuer}
\end{align}
Defining a difference in positions per unit of time as
\[
\varepsilon(n)=n(q(n)-p(n+1))\,,
\]
we have using \eqref{E:pursuer} at time $t_{n+2}$, and
\eqref{E:target} at time $t_{n+1}$ that:
\begin{equation}
\varepsilon(n+1)=\varepsilon(n)+\gamma(n+1)-V(n+2)\,.
\end{equation}
The reader might well anticipate that the pursuit strategy is the
greedy algorithm: namely, the predator chooses vertex $V(n+1)$
which minimizes of $\varepsilon(n)$ with respect to the Euclidean
norm, with some tie breaking rule. In terms of the notation
previously introduced
\[
V(n+1)=v(\varepsilon(n-1)+\gamma(n))=v((n-1)(q(n-1)-p(n))+\gamma(n))\,.
\]
Thus
\begin{equation}
\varepsilon(n+1)=\varepsilon(n)+\gamma(n+1)-v(\varepsilon(n)+\gamma(n+1))\,,
\end{equation}
which we recognize as Equation \ref{recursion}. It then follows
from Theorem \ref{T:main} that
\[
||q(n)-p(n+1)||\rightarrow 0\,,
\]
which implies
\[
||q(n)-p(n)||\rightarrow 0\,.
\]
In this sense the predator catches the prey.

\subsection{Instances of the greedy algorithm in classical mathematics.}
The greedy algorithm with which we are concerned in this paper
leads one into a realm of tremendous richness. Consider the
simplest input case where $\gamma(k)\equiv \gamma$ is constant.
The performance of the greedy algorithm reduces to the study of
iterations of a single map $\phi_\gamma.$ This map is a piecewise
isometry.  In general the analysis of piecewise isometries can be
of incredible depth and difficulty (see for example \cite{AdKiT}).
Indeed, the simplest of all cases, the constant input case where
the polytope $P=[0,1]$ the unit interval has significant
mathematical substance. One easily checks that the interval
$I=[\gamma - \frac{1}{2}, \gamma + \frac{1}{2}]$ is invariant and
absorbing under $\phi_\gamma$. Both the interval $I$ and $[0,1]$
are fundamental regions for the action of $\mathbb{Z}$ on
$\mathbb{R}$ and
$$\phi_\gamma(x)=x+\gamma-v(x)= x+\gamma \mod{\mathbb{Z}}.$$
In other words, both $I$ and $[0,1]$ can be identified with the
circle $\mathbb{R}/\mathbb{Z}$ and $\phi_\gamma$ with a rotation
by an angle $\gamma$. This is one of the standard examples studied
in ergodic theory and has the property of \textit{unique
ergodicity} when $\gamma$ is irrational.
\begin{figure}[ht]
 \centerline{\epsfbox{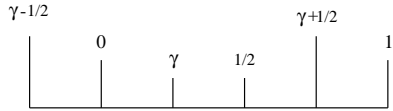}}
 \caption{Diagram for simplest case}
 \label{intervals}
\end{figure}
In Figure 1 we have drawn the case in which $0<\gamma<1/2.$  A
consequence of unique ergodicity is that for irrational $\gamma$
and any $x_0\in I$, the ergodic averages of the characteristic
function $\chi_{[\frac{1}{2},\gamma
+\frac{1}{2}]}(x_0+k\gamma\mod{1})$ converges to the measure of
$[\frac{1}{2},\gamma +\frac{1}{2}]$ which of course is equal to
$\gamma$. For this particular interval, \eqref{E:averages} and the
equation
$$\chi_{[\frac{1}{2},\gamma +\frac{1}{2}]}(x_0+k\gamma\mod{1})=
\chi_{[\frac{1}{2},\gamma +\frac{1}{2}]}(\phi^k_\gamma (x_0)=
v(x(k)),$$ where $x(0)=x_0,$ leads to the same result even for
rational $\gamma.$

\medskip
The type of sequence $v(x(k))$  of 0's and 1's generated by the
recursion \eqref{recursion2} is known as a \textit{Sturmian
sequence.} Such sequences were studied by Morse and Hedlund
\cite{MH}. They have the property that there is a limiting
frequency of 1's and the distribution of 0's and 1's is as even as
possible. Arbitrary Sturmian sequences can be identified with
orbits of homeomorphisms of the circle, but only those that appear
for rigid rotations get generated by the greedy algorithm in the
constant input case (see \cite{MH} and \cite{SieTZ}). These
sequences provide a symbolic counterpart to the theory of
Diophantine approximation (see \cite{SieTZ} and \cite{LagTre}) and
have a long history. Johan Bernoulli was apparently the first to
discover what came to be known as Sturmian sequences. He did that
to describe an easier method of interpolation in astronomical
tables\cite{Bern}; Christoffel applied such sequences to modular
arithmetic \cite{Christ}; Smith used them to study well
distributed sequences of two symbols \cite{Smith}. Both
Christoffel and Smith seemed to have been unaware of the earlier
work of Bernoulli. Bernoulli's work was later revisited by Markov
\cite{Mark1882}. A generalization of Smith's view of the problem
was considered in \cite{BRTT} where the circle was replaced by
arbitrary dimension tori, and the corresponding multidimensional
continued fractions were described.

\medskip
One more feature of the constant input case is the fact that under
the action of $\phi_\gamma$ on $\mathbb{R}$ the interval I is
invariant and \textit{absorbing}--\textit{i.e.,} every orbit
eventually enters and stays in I. We expect that the following
generalization holds in all dimensions, the proof of which we
shall defer to a later work. For any lattice $\mathcal{L}$ with a
preferred basis in $\mathbb{R}^N,$ let $P$ be the standard
fundamental parallelepiped in the basis. Then for any $\gamma \in
P$ there exists another fundamental region $Q_\gamma$ which is an
invariant absorbing set under the action of $\phi_\gamma $; and on
$Q_\gamma$ the map $\phi_\gamma$ is a rotation by $\gamma$ on
$\mathbb{R}^N/\mathcal{L}.$ The same results holds if $P$ is a
simplex instead of a parallelepiped leading however to a different
fundamental region $Q_\gamma.$ Furthermore, there is a convergence
of ergodic averages for special sets in the simplex case similar
to the interval case. Finally, it is conceivable that this result
could lead to higher dimensional analogues to Sturmian sequences.

\subsection{Outline of the paper.}
In \S2 we describe how the proof of Theorem \ref{T:main} evolved
from dimension 2 to the general case. In \S3 we prove Theorem
\ref{T:main} based on two propositions, Proposition
\ref{proposition1} and Proposition \ref{proposition2}, which we
prove in \S4. In \S4, we show that the proof of Proposition
\ref{proposition2} itself easily follows from a third proposition,
Proposition \ref{proposition3}. Proposition \ref{proposition3}
posits the existence of a special subset of the
$(N-1)$-dimensional sphere upon which the construction of the
invariant set $Q$ is based. We postponed the proof of that
Proposition to \S5. Item \textit{(3)} of that proposition is where
all the difficulties and subtleties of our proof of Theorem
\ref{T:main} lie. Finally, at the end of \S5 we have finished
assembling all the elements for a complete proof of Theorem
\ref{T:main}, and in \S6 we present two Theorems concerning some
additional properties which are byproducts of our method of proof.

% Section 2
% Section 2
\section{How the proof evolved from dimension 2 to the general case}
% Section 2
% Section 2

We shall adhere to the convention of writing the inner product of
vectors as a dot product, the Euclidean norm as $|| \cdot ||,$ and
Euclidean distance between a pair of points or between a point and
a set as $d(\cdot,\cdot)\, .$

\medskip

\begin{definition}\label{D:Vor} Given a set of points  $v_0,\dots,v_{M-1}$
in any dimension, the \textit{Vorono\"{\i} region} $R_{v_i}$ of
$v_i$ is defined as
\[
R_{v_i} = \bigcap_j\{ x: ||x-v_i|| \leq ||x-v_j|| \},
\]
hence a polyhedron being the finite intersection of half-spaces.
It is the set of points closer to $v_i$ than to any other point in
the set (or possibly no further than equidistant to some). Also
let $v(x)$ denote\textit{ the closest vertex to $x$ with a tie
breaking rule }such as choosing the vertex of smallest index in
case of ties.
\end{definition}

Now let $P$ be a convex polygon with vertices $v_0,\dots,v_{M-1}\,
.$ To make the definition of the map $\phi_\gamma(x)$ more
specific, the tie breaking rule for $v(x)$ will be for us the
smallest index $i$ such that $x\in R_{v_i}.$ For example, if $P$
is an interval $[v_0,v_1],$ then the Vorono\"{\i} regions are:
\begin{align}
R_{v_0}&=(-\infty, \frac{v_0+v_1}{2}]\,,\notag  \\
\intertext{and} R_{v_1}&=(\frac{v_0+v_1}{2},\infty )\,.\notag
\end{align}
The mapping $\phi_{\gamma}:\mathbb{R}\hookrightarrow\mathbb{R}$ is
given by:
\begin{equation}\label{rot} \phi_{\gamma}(x)=x+\gamma-v(x), \text{
where } v(x)=
\begin{cases}
v_0, \text{ if } x\leq \frac{v_0+v_1}{2}\,,\\
v_1, \text{ if } x > \frac{v_0+v_1}{2}\,.\notag
\end{cases}
\end{equation}

\medskip
Theorem \ref{interval} and Theorem \ref{polygon} are respectively
dimension 1 and 2 (intervals and polygons) versions of Theorem
\eqref{T:main}. Proofs of them appeared in \cite{AKMTW2002}.

\begin{theorem}\label{interval}
Let $Q_t$ denote the interval generated by moving the ends of
$P=[v_0,v_1]$ outward a distance $t\geq 0:$ \textit{i.e.},
$Q_t=[v_0-t,v_1+t].$ Then for any $t\geq \frac{v_1-v_0}{2}$ and
$\gamma \in P$
\[
\phi_{\gamma}(Q_t)\subset Q_t.
\]
\end{theorem}

\bigskip Let  $v_i,\ i=0,\dots, M-1,\ n \geq 2$ indexed in
clockwise order $P$ and $n_j,\ j=1,\dots,n$ the unit normal
vectors to edges $v_{j-1}v_{j}$ of $P$. Then $n_j\cdot
(v_j-v_{j-1})=0,$ and we can write $n_j\cdot v_{j-1}=n_j\cdot
v_{j}=d_j.$

\begin{figure}[ht]
 \centerline{\epsfbox{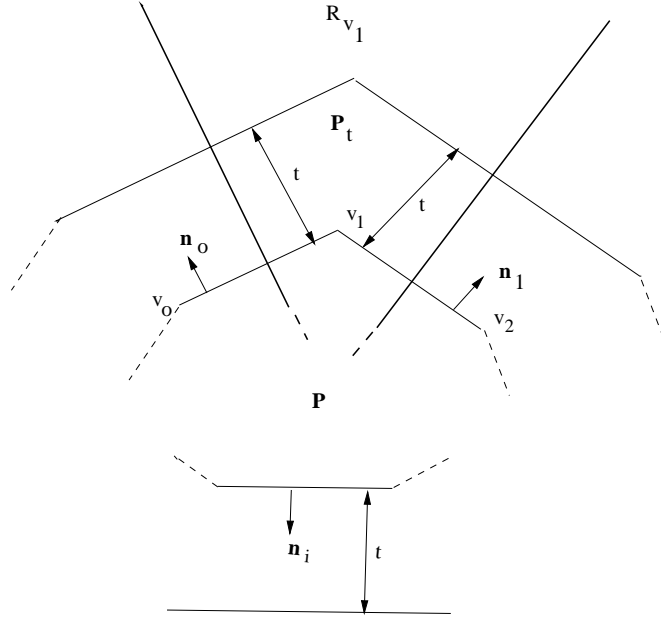}}
 \caption{Polygons $P,\ P_t,$ and Vorono\"{\i} region $R_{v_{1}}$}
 \label{F:polygons}
\end{figure}

\begin{theorem}\label{polygon}
Let $Q_t$ denote the polygon generated by moving the edges of $P=
\bigcap_j \{ x:x\cdot n_j \leq d_j\}$ perpendicularly outward a
distance $t\geq 0:$ \textit{i.e.},
\[
Q_t = \bigcap_j \{ x:x\cdot n_j \leq d_j + t\}.
\]
See Figure \ref{F:polygons}. There exists $T$ such that for any $t
\geq T$ and $\gamma \in P$
\[
\phi_{\gamma}(Q_t)\subset Q_t.
\]
\end{theorem}

\bigskip
The above method of constructing an invariant $Q$ by outwardly
translating the planes of faces of a polytope $P,$  even by
different distances, is doomed to failure in dimension $N=3.$ as
shown by the following counterexample.

\begin{figure}[ht]
 \centerline{\epsfbox{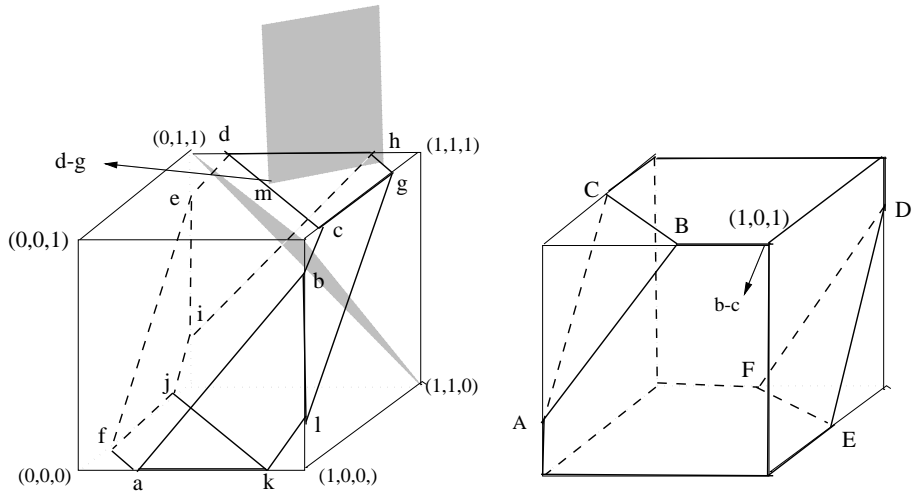}}
 \caption{Polytope Counterexample}
 \label{cube}
\end{figure}

Consider the octahedral polytope $abcdefghijkl$ in Figure
\eqref{cube}. The faces $abcdef$ and $ghijkl$ are parallel to the
plane $(0,0,0)(0,1,1)(1,0,1)$ cutting inside the cube as shown and
equidistant from the points $(0,0,0)$ and $(1,1,1)$ respectively.
The face of the Vorono\"{\i} regions $R_g$ separating $g$ from $h$
is in the plane perpendicular to the top face of the cube through
the midpoints of segments $dc$ and $hg.$  Invariance fails at the
point the midpoint $m$ of $dc$ because here the vector $d-g$
sticks out of the polytope. No matter how far the top plane of the
cube is translated upward this vector will still stick out the
resulting polytope at such a point. So this means that the faces
$abcdef$ and $ghijkl$ must be translated far enough apart. But
then a new difficulty appears: the face of Vorono\"{\i} region
$R_c$ lies in the plane determined by $(0,1,1)$, $(1,0,1)$, and
the midpoint of $bc.$ The vertex $(1,0,1)$ of the cube lies in
this face, but here the vector $b-c$ sticks out of this new
candidate for an invariant polytope. On one hand this second
difficulty cannot be overcome by moving the front and side faces
outward, but on the other hand moving the top face up reintroduces
the first difficulty. We leave it to the reader to verify that one
problem cannot be overcome without introducing the other: for
details and the general case of dimension $N>2$, see
\cite{GAGA:IIUnavoidableDifficulties:2003}.

\medskip
So another idea for a general construction of an invariant set is
needed (special constructions for special polytopes will be
discussed in \cite{GAGA:IVConstInvReg:2003}). Imagine the limiting
case of a sphere centered at the origin scaled down to the unit
sphere as the radius goes to infinity. The edges at a vertex are
normals to hyperplanes through the origin, which determine the
limiting Vorono\"{\i} region at that vertex. The following, which
we shall verify later for all dimensions, is easy to see in the
two dimensional case below.  It is the most important element of
the proof. \textit{The cone of normals to the supporting
hyperplanes at a vertex agrees with the Vorono\"{\i} region at
infinity of that vertex} (See proof of Proposition
\ref{proposition1}). Thus the tangent space of the sphere at a
normal in the Vorono\"{\i} region of a vertex is a supporting
hyperplane of the polytope at that vertex. In this idealized
situation the polytope has shrunk to a point at the origin; and
the unit sphere (equivalently the sphere at infinity) is
invariant. Backing off from infinity to a sufficiently large
sphere, the invariance fails only on neighborhoods of Vorono\"{\i}
boundaries. Resolving this difficulty then becomes the main task
of the proof. To do this we construct a special subset $\Omega$ of
the unit sphere the tangent planes of which are supporting
hyperplanes of a \textit{proto-invariant} convex set $Q_\infty.$
Then the set $Q=\rho Q_\infty$ for large enough $\rho$ will be the
desired invariant convex set.

\bigskip
 Before launching into the formal proof for all dimensions, we
illustrate the idea for some low dimensional polytopes indicating
the induction upon which the general construction is based.

\begin{figure}[ht]
 \centerline{\epsfbox{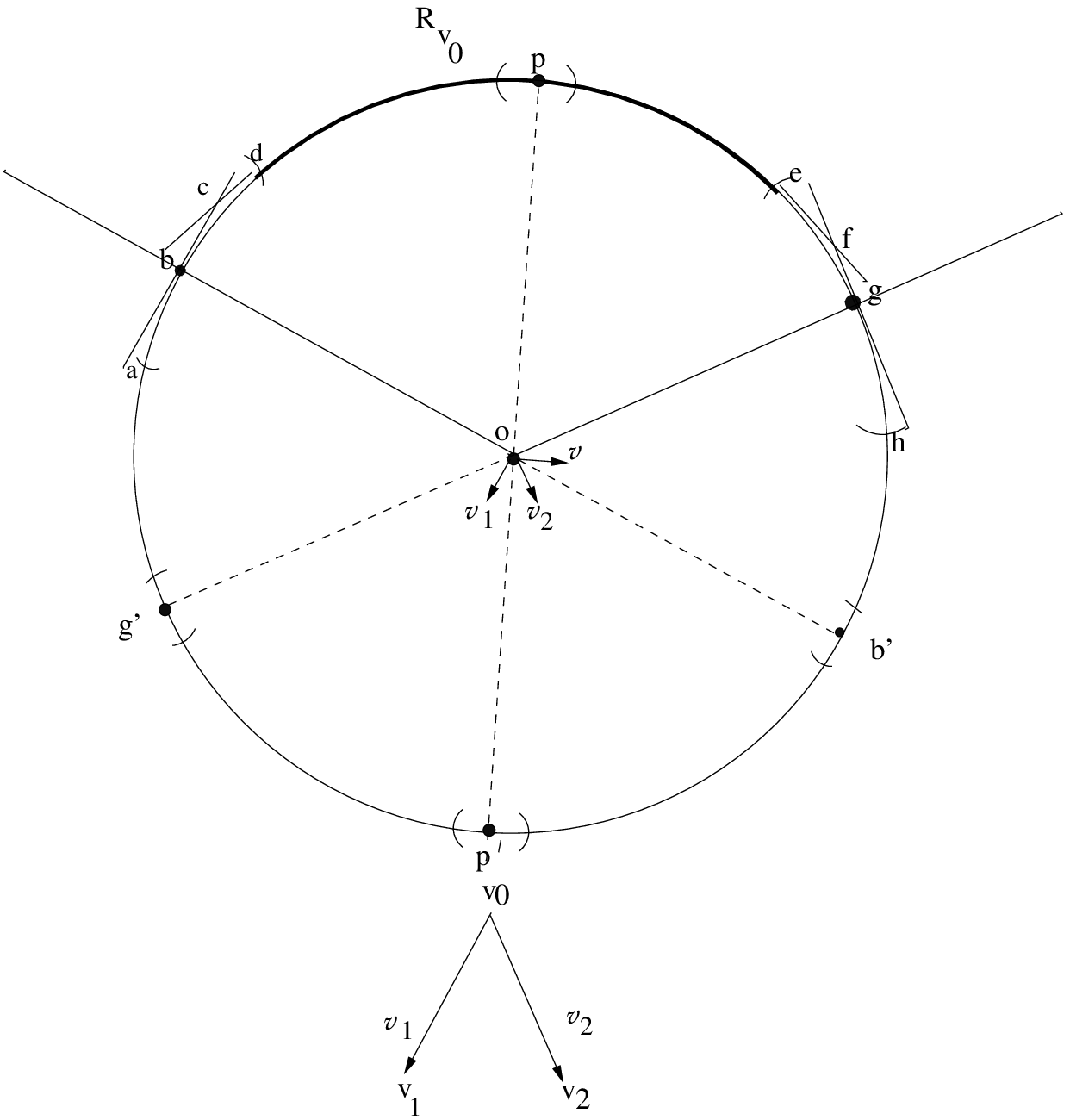}}
 \caption{Construction of $Q$ in dimension2}
 \label{circle}
\end{figure}

First we take up a construction of Of $Q$ viala set$\Omega$ in
dimension 2. Consider the idealized limiting circle scaled down to
the unit circle as depicted in Figure \eqref{circle}.
Let $v_1,v_0,v_2$ be three successive vertices of a polygon, and
$\nu_1=v_1-v_0$ and $\nu_2=v_2-v_0$. The points $b$ and $g$ lie
the intersection of rays perpendicular to $\nu_1$ and $\nu_2$ with
$S.$ These rays are parallel to the unbounded edges of the of the
Vorono\"{\i} region $R_{v_0}$ for the vertex $v_0.$ Using the same
letters to denote corresponding points on a dilated circle $\rho
S$, the intersection of the infinite bounding edges of the actual
Vorono\"{\i} region with the circumference of the dilated circle
for large $\rho$ would be in the arc $ad$ and in the arc $eh.$ The
only place a vector $\gamma-v_0$ would stick out of the circle for
some $\gamma \in P$  would be at points $b$ and $g$
(\textit{e.g.}, $\gamma -v_1$ at $b$) and nearby points relatively
speaking. So in order to avoid this difficulty, we remove small
open arcs $ad$ and $eh$ of $S,$ retaining the points $b$ and $g,$
and carry out this procedure for each Vorono\"{\i} region
$R_{v_j}.$ The size of the arcs is subject to the restriction that
the removed arcs don't interfere with one another. It can then be
seen that, if $Q_\infty$ is the convex set supported by tangent
lines to the unit sphere at points of $\Omega$, then the convex
set $\rho Q_\infty$ is invariant for large enough $\rho.$

\begin{figure}[ht]
 \centerline{\epsfbox{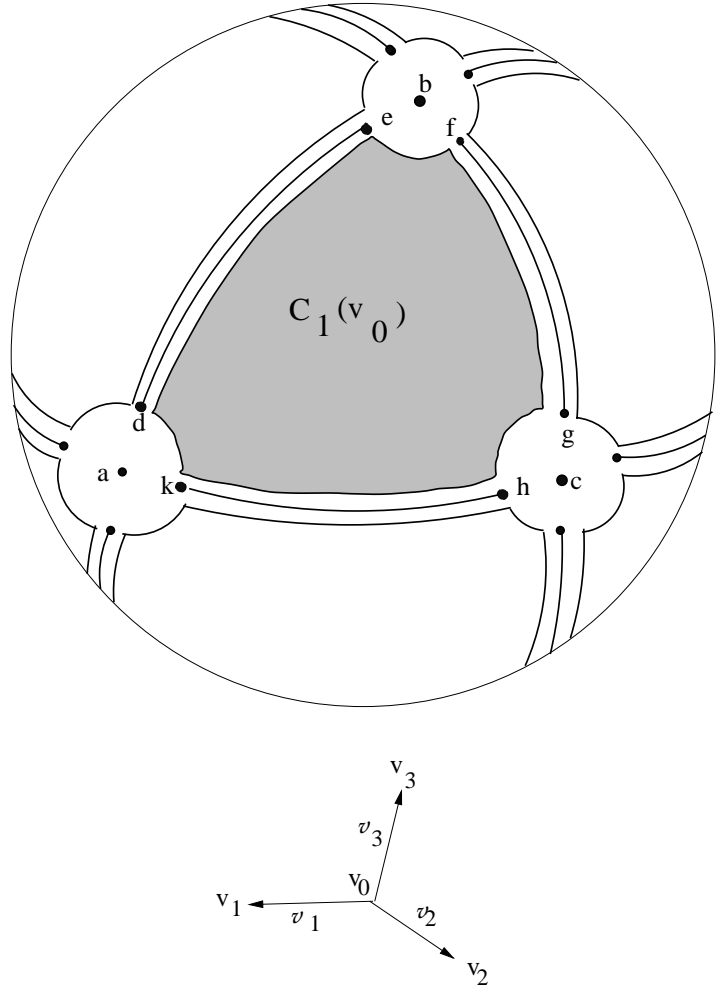}}
 \caption{Construction of $Q$ in dimension 3}
 \label{sphere}
\end{figure}

\medskip
A construction for $Q$ via  a set $\Omega$in dimension 3 is
similar. For the sake of simplicity, suppose that three edges meet
at vertex $v_0$ and that $v_1,$ $v_2,$ and $v_3$ are the
neighboring vertices, as in the bottom of Figure
\eqref{sphere}.\textsf{}
Let $\nu_1=v_1-v_0,$ $\nu_2=v_2-v_0,$ and $\nu_3=v_3-v_0$ be the
normals to the faces of the the Vorono\"{\i} regions $R_{v_0}.$
The planes defined by these normals intersect the idealized scaled
limiting unit sphere in the spherical triangle $abc.$ To construct
$\Omega$ we remove from the unit sphere small open spherical caps
centered at $a,$ $b,$ and $c,$ but retain these centers. In the
complement of what has been removed, we remove open strips
centered on the edges of the spherical triangle, but retain their
center lines which in this case are the closed arcs $de,$ $fg,$
and $hk.$ In Figure \eqref{sphere} we have illustrated the
construction for Vorono\"{\i} region $R_{v_0}.$ We carry out this
procedure for each Vorono\"{\i} region. Again caps and strips
should be chosen so that they do not interfere with one another.
The set $\Omega$ is what survives. As before, it can be seen that
the failure of invariance of large spheres is overcome by the set
$\Omega,$ and the convex set $Q_\infty$ supported by tangent
planes to the unit sphere at points of $\Omega,$ whereupon $Q
=rQ_\infty$ is invariant for large enough $r.$

\medskip
Note that the boundary of the removed caps from the 2-sphere in 3
dimensions has the same structure as $\Omega$ in 2 dimensions.
Also the boundary of a cross section of a removed strip from the
sphere has the same structure as a line through a boundary of a
removed arc on the circle (like the situation on the arc $ad$).
Thus we get a glimpse of the induction about to be carried out.

\medskip
In dimension 4, the intersection of a Vorono\"{\i} region and the
3-sphere is a spherical polytope. To construct $\Omega$ we proceed
exactly as before. We remove small open 4-spheres about the
vertices of the polytope retaining the centers. From what remains,
we remove open tubes around the edges, retaining the edges.
Finally from what still remains of the unit sphere we remove slabs
about the faces of the polytope while retaining the faces.  We do
this for each Vorono\"{\i} region. Again a restriction on the
removed sets is that they don't interfere with one another.

\medskip
Again we note that the boundary of a removed spherical
neighborhood from the 3-sphere in 4 dimension is a 2-sphere having
the same structure as the 3-dimensional $\Omega.$ The retained
center of such a removed spherical neighborhood is not visible in
the 3-dimensional subspace in which the bounding 2-sphere lives.
Rather the center of this 2-sphere is a projection of the retained
point in $\Omega.$ This corresponds to what is happening in the
previous dimension.  In Figure \eqref{sphere} the point $a$ is not
in the plane which contains the boundary circle of the removed
spherical cap on the 2-sphere. The center of this circle is the
projection of the point $a$ to this plane. Going down one more
dimension, this corresponds to the fact that $b$ does not lie on
the line through $a$ and $d.$  A cross section of a removed tube
has the same structure as the boundary of a removed cap in
dimensional 3 which in turn has the same structure as the removed
circle in dimension 2.

\medskip
\textbf{Remark.} Returning to Figure \eqref{circle} we note more
points could have been removed between $d$ and $e$ without
compromising invariance. Similarly, the same can be said for
removing more points from $\Omega$ in any higher dimension. Our
actual construction of $\Omega$ will take advantage of this fact,
and we shall remove much more from the unit sphere than what we
have just described. This has two virtues. First describing
removed sets is slightly easier. Second it leads to a simplified
proof from which we get for free Theorem \ref{T:byproduct2}.

\medskip
For the construction in dimension 2 of this more truncated
$\Omega$, let $p$ and $p^\prime$ be the points of intersection of
the unit sphere with the line normal to $\nu={\rm v}_i-{\rm v}_j$
for a pair of vertices ${\rm v_i}$ and ${\rm v}_j.$ Remove from
the unit sphere small open arcs (symmetric) about $p$ and
$p^\prime$ but retaining $p$ and $p^\prime.$ Do this for
\textit{every} pair of vertices ${\rm v}_i$ and ${\rm v}_j$ of the
polytope. The only restriction is that removed arcs be small
enough not to interfere with each other. The set $\Omega$ is
smaller than in the previous construction which results in fewer
supporting lines for $Q_\infty$ which then results in a larger
$Q_\infty$. Nevertheless, the ultimate $Q$ will still be
invariant.

\medskip
To construct $\Omega$ in dimension 3, consider all possible
intersections with the unit sphere of planes whose normals are
given by ${\rm v}_i-{\rm v}_j$ for every distinct pair of adjacent
vertices ${ \rm v}_i \ne {\rm v}_j$ of the polytope $P$. These
intersections consist of one dimensional lines, and two
dimensional planes.  We first remove caps centered at all the
points that are intersections of the one dimensional lines with
the unit sphere while retaining the center points. From what is
left we remove all strips centered at great circle arcs that are
intersections of two dimensional planes with the unit sphere while
retaining the central arcs. Once again the removed sets must be at
least small enough so that they do not interfere with one another.
To facilitate our proof it turns out to be convenient to make
these removed sets even smaller than just the requirements for non
interference. As before, the resulting $\Omega$ leads to an
invariant $Q.$

\medskip
The above will be the procedure for constructing $\Omega$ in all
dimensions.  Showing that a supporting hyperplane tangent to the
unit sphere at a point in $\Omega$ is not cut off at the wrong
place by another one is delicate. This is clear in dimension 2,
but far from obvious in higher dimensions.

% Section 3
% Section 3
\section{Proof of Theorem \eqref{T:main}}
% Section 3
% Section 3
We shall use freely some notations, definitions, and basic results
from classical convex geometry (see for instance
\cite{rockafellar:convex_analysis} or \cite{S} and references
therein).

\subsection{Notation}
\begin{itemize}
     \item $P$ is a polytope--namely, the convex hull of a finite set
     (vertices) $\{ v_1, \dots,v_M \}$ of extreme points
     (vertices)in $\mathbb{R}^N$. By the classical \textit{Finite basis theorem for polytopes}
     attributed to Minkowski [1896], Steinitz [1916], and Weyl [1935]
     (see \cite{S}, p. 89), $P$ is a bounded polyhedron.
     \item $C(v_i)$ denotes the cone of outward normals to the supporting
       hyperplanes of $P$ at vertex $v_i;$ and $C_1(v_i)$ the
       unit normals in $C(v_i).$ See spherical triangle $abc$ in
       Figure\eqref{sphere}.
    \item As a polyhedron, $P$ can be expressed by
      \begin{equation}\label{E:P}
      P=\bigcap_{\text{all }v_i}\bigcap_{\omega\in C(v_i)}
      \{ x| \omega \cdot (x-v_i) \le 0\}\,.
     \end{equation}

\end{itemize}

\bigskip
\subsection{Proof of the main result}
We shall assume two propositions, the proofs of which shall be
deferred to the next two sections.

\medskip
\begin{proposition}\label{proposition1}
For each $\varepsilon >0$ there exists $D>0$ such that if $||x||>
D$ then $d(\frac{x}{||x||},C_1(v(x)))<\varepsilon.$
\end{proposition}

\bigskip
\begin{proposition}\label{proposition2}
There exists a closed subset $\Omega$ of the unit sphere in
$\mathbb{R}^N$ centered at the origin with the following
properties:
\begin{enumerate}
\item $\Omega$ is not confined to any half-sphere;
\item there exists an $\varepsilon,\ 0<\varepsilon<1, $
such that if $y\in C_1(v)$ and $\omega \in \Omega -C_1(v) $ then
there exists $\omega^\prime \in C_1(v)\cap \Omega$ such that
\begin{equation}\label{E:inequality}
d(\omega^\prime, y)<d(\omega,y)-\varepsilon\,.
\end{equation}
\end{enumerate}
\end{proposition}

\bigskip
\noindent \textbf{Proof of Theorem \ref{T:main}.}

\bigskip
 Let
\begin{align}
Q_\infty &\equiv \bigcap_{\omega\in \Omega}\{ x| x\cdot \omega \le 1 \}\,, \\
\intertext{and} \rho Q_\infty &=\bigcap_{\omega\in \Omega}\{ x|
x\cdot \omega \le \rho \}\,.\label{E:rhoQinfty}
\end{align}
We shall show that for all sufficiently large $\rho$ the set $\rho
Q_\infty$ satisfies Theorem \ref{T:main}. It is clear by
definition that $Q_\infty$ and hence $\rho Q_\infty$ are convex.
If $Q_\infty$ were not bounded, it would contain a ray
$\mathbb{R}^+x=\{ tx|\ t\ge 0 \},$ and the half-sphere determined
by the hyper-plane normal to $x$ would contain $\Omega,$
contradicting Proposition \ref{proposition2}. Thus $Q_\infty$ and
$\rho Q_\infty$ are bounded.

\medskip
We must show that there exists an \textsf{R} such that
\begin{equation}\label{E:Q}
\omega \cdot\phi_\gamma(x)=\omega\cdot(x+\gamma-v(x)\le \rho
\end{equation}
for all $x\in \rho Q_\infty,\ \rho \ge \textsf{R},\ \gamma \in P,
\text{ and } \omega\in \Omega.$

\bigskip
\noindent Case I: $\omega \in C(v(x)).$

\medskip
By \eqref{E:P} and \eqref{E:rhoQinfty}, for $\gamma \in P$ we have
\begin{equation}\label{E:inequality2}
\omega \cdot x+\omega \cdot (\gamma-v(x)) \le \rho\,.
\end{equation}
The inequality is satisfied no matter what the size of $||x||$ and
$\rho>0$ are.

\bigskip
\noindent Case II: $\omega \notin C(v(x)).$

\medskip
Let $\varepsilon$ be the one given in Proposition
\ref{proposition2}. From Proposition \ref{proposition1} there
exists  $\textsf{R}_0$ such that for $||x||>\textsf{R}_0$ there
exists $y \in C_1(v(x))$ satisfying
\[
d(\frac{x}{||x||},y))< \varepsilon/4\,,
\]
By Proposition \ref{proposition2} there exists $\omega^\prime \in
C_1(v(x))\cap \Omega$ such that
\[
d(w,y)>d(\omega^\prime,y)+\varepsilon\,.
\]
By the triangle inequality
\begin{align*}
d(\omega, \frac{x}{||x||})&\ge d(\omega,y)-d(y,\frac{x}{||x||})\\
&> d(\omega^\prime,y)+\varepsilon -d(y,\frac{x}{||x||})\\
&\ge d(\omega^\prime,\frac{x}{||x||})+\varepsilon -2 d(y,\frac{x}{||x||})\\
&\ge d(\frac{x}{||x||},\omega^\prime)+\varepsilon/2 \,.
\end{align*}

\noindent Let $\theta$ be the angle between $x$ and $\omega,$ and
$\theta'$ between $x$ and $\omega'.$ Using the fact that
differences in arc lengths are greater than differences in
corresponding chord lengths (an exercise--hint: make chords
parallel) we get
\begin{equation}
\theta-\theta' > d(\omega, \frac{x}{||x||})-
d(\frac{x}{||x||},\omega^\prime) > \varepsilon/2
\end{equation}
for the arc length difference. Thus
\begin{equation}
\cos(\theta)<\cos(\theta'+\varepsilon/2)<\cos(\theta')\cos(\varepsilon/2)\,,
\end{equation}
and we arrive at
\[
\omega \cdot \frac{x}{||x||}< \omega^\prime \cdot
\frac{x}{||x||}(1-e)
\]
where $e=1-\cos{\frac{\varepsilon}{2}}.$ Note $\ 0<e<1,$ depending
only on $\varepsilon.$ We then get
\begin{equation}\label{E:inequality3}
\omega \cdot x+\omega \cdot (\gamma-v(x)) <\omega^\prime\cdot
x(1-e)+\text{diam}(P)\,.
\end{equation}
We can therefore choose $\textsf{R}_1>\textsf{R}_0$ such that for
$||x|| \ge \textsf{R}_1$ the right hand side of
\eqref{E:inequality3} is $< \rho.$ Let
\begin{equation}\label{E:R}
\textsf{R}=\textsf{R}_1+\text{diam}(P)\,,
\end{equation}
and let $x\in \rho Q_\infty$ where $\rho >\textsf{R}.$ By what we
have just shown, if $||x||\ge \textsf{R}_1$ then $x$ satisfies
equation \eqref{E:Q}: so $\phi_\gamma(x)\in \rho Q_\infty.$
Finally, if $||x||<\textsf{R}_1$ then $x\in \rho Q_\infty:$
because $||x-\phi_\gamma(x)||\le \text{diam}(P).$ $\Box$

\medskip
\begin{numberedremark}
\bigskip At this point we have also proved that given any compact set $K$ the set $Q$ can be chosen to contain
to contain $K$.
\end{numberedremark}

% Section 4
% Section 4
\section{Proofs of Propositions \eqref{proposition1} and \eqref{proposition2}}
% Section 4
% Section 4

 \bigskip
\subsection{Proof of Proposition \ref{proposition1}}

By a standard argument in linear programming--namely, for each
vector $u$ the function $u\cdot x $ achieves a maximum over $P$ at
some vertex $v_i$ which means that $u\cdot (x-v_i)\le 0$ for $x\in
P,$ or equivalently that $u\in C(v_i)$-- we have $\mathbb{R}^N =
\bigcup_i C(v_i).$

Furthermore, by the  \textit{Decomposition Theorem of Polyhedra}
due to Motzkin [1936] (see \cite{S}, p.88) a Vorono\"{\i} region
$R_{v_i}$, being a polyhedron, can be written
\[
R_{v_i}=P_{v_i} + K_{v_i}\,,
\]
where $P_{v_i}$ is a polytope (recall a polytope is a bounded set)
and $K_{v_i}$ is a polyhedral cone. Since $\mathbb{R}^N =
\bigcup_i R_{v_i}$, it is easy to prove that $\mathbb{R}^N =
\bigcup_i K_{v_i}$ as well.

We show that $K_{v_i}=C(v_i)$ as follows. Let $u\in C(v_i),$ and
$t\ge 0.$ The vector $u$ determines a supporting hyperplane for
$P$ at $v_i$ which implies that $v_i$ is the closest vector of $P$
to the vector $tu+v_i:$ so $tu+v_i \in R_{v_i}.$ Thus from the
polyhedral decomposition theorem $tu+v_i=p(t)+k(t)$ for all $t\ge
0$ where $p(t) \in P_{v_i}$ and $k(t)\in K_{v_i}.$ Therefore,
\begin{equation}\label{E:cone}
u=\frac{p(t)-v_i}{t} +\frac{k(t)}{t}\,.
\end{equation}
The first term on the left in equation \eqref{E:cone} converges to
0 as $t\rightarrow \infty$ so that the second term, which always
belongs to $K(v_i)$ by convexity, is bounded for all $t>0.$ Thus
there is a convergent subsequence whose limit is in $K(v_i)$ which
implies that $u\in K_{v_i}:$ hence $C(v_i)\subset K_{v_i}.$

Since the two collections of $C(v_i)$'s and of $K_{v_i}$'s form
partitions of $\mathbb{R}^N$, the reverse inclusion follows easily
from the fact that the interiors of $C(v_i)$ and $C(v_j)$ are
disjoint for $i\ne j,$ and similarly for the interiors of
$K_{v_i}$ and $K_{v_j}$.

\medskip Given any $x \in \mathbb{R}^N,$ we have just shown
\[
x\in P_{v(x)}+C(v(x))\,,
\]
so there exists $p(x)\in P_{v(x)}$ such that $x-p(x)\in C(v(x))$.
Theorem \ref{proposition1} now follows from the easily proved
limit
\[
||\frac{x}{||x||}-\frac{x-p(x)}{||x-p(x)||}|| \rightarrow 0 \ ,
\text{as } ||x|| \rightarrow \infty\,. \ \Box
\]

\bigskip
\subsection{Proof of Proposition \ref{proposition2}}

\medskip
The proof of Proposition \ref{proposition2} follows easily from
the following proposition, the proof of which has been postponed
until the last section. As previously stated, this proposition is
the heart of the matter

\medskip
\begin{proposition}\label{proposition3}
There exist a partition of each $C_1(v_i)$ into a finite number of
connected sets $\sigma$ and  a closed subset $\Omega$ of the unit
sphere $S$ in $\mathbb{R}^N$ centered at the origin with the
following properties:
\begin{enumerate}
    \item $\Omega$ is not confined to any half-sphere;
    \item for each $\sigma \subset C_1(v_i)$ there is a
    $\epsilon_\sigma>0$ such that
    \[\Omega\cap {\rm U}_{\varepsilon_\sigma}(\bar{\sigma}) \subset
    \bar{\sigma},
    \]
    where  ${\rm U}_\epsilon(A)=\{x| \ d(x,A)< \epsilon \}$ denotes
   the $\epsilon-$neighborhood of a set $A;$
    \item the closest points $\omega \in \Omega$ to $y \in \bar{\sigma}
    \subset C_1(v_i)$ are in the same set closure $\bar{\sigma}$.
    In particular $\bar{\sigma} \cap \Omega \ne \emptyset.$
\end{enumerate}
\end{proposition}

The sets $\sigma\in C_1(v_i)$ will be chosen as \textit{cells} of
a simplicial decomposition induced by great spheres on $S$ as
specified in \S5.

\bigskip
\noindent \textbf{Proof of Proposition \ref{proposition2}.}

\medskip
We must prove that there exists $\varepsilon$ such that for all
$v_i$ and all $y\in C_1(v_i)$
\begin{equation}\label{E:epsilon} d(y,\Omega \setminus
C_1(v_i))-d(y,C_1(v_i)\cap \Omega)>\varepsilon\, .
\end{equation}

\medskip
For $\sigma \subset C_1(v_i)$ define a function
$f_\sigma:\bar{\sigma}\to \mathbb{R}$ by
\[
f_\sigma (y) \equiv d(y,\Omega \setminus {\rm
 U}_{\epsilon_\sigma}(\bar{\sigma}))-d(y,C_1(v_i)\cap \Omega)
\]
where $\epsilon_\sigma$ is given by Proposition
\ref{proposition3}.  Then $f_\sigma$ is continuous and strictly
positive on a compact set: hence assumes a minimum $e_\sigma>0.$
We claim $\varepsilon=\min e_\sigma$, where the minimum is taken
over all $\sigma\subset C_1(v_i)$ and all $v_i$'s, satisfies
\eqref{E:epsilon}. Observe that
\begin{equation}\label{E:inequality4}
d(y,\Omega \setminus C_1(v_i))-d(y,C_1(v_i)\cap \Omega) \ge
d(y,\Omega\setminus \bar{\sigma})-d(y,C_1(v_i)\cap \Omega)\, ,
\end{equation}
where $y\in \bar{\sigma}\subset C_1(v_i).$ Now by (2) and (3) of
Proposition \ref{proposition3} the right hand side of
\eqref{E:inequality4} is equal to $d(y,\Omega \setminus {\rm
U}_{\epsilon_\sigma}(\bar{\sigma}))-d(y,C_1(v_i)\cap
\Omega)=f_\sigma(y)\ge \varepsilon.$ $\Box$

% Section 5
% Section 5
\section{Proof of Proposition \ref{proposition3} }
% Section 5
% Section 5

\subsection{Notation}
\begin{itemize}
     \item $S=S(c,r)\subset \mathbb{R}^{n+1} = \mathbb{R}^N$ denotes an
     $n$-dimensional sphere
     with radius $r$ and center $c$. We shall call the intersection of $S$ with any
     $k$-dimensional affine subspace of $\mathbb{R}^{n+1}$ containing $c$ a
     \textit{great sphere of dimension $k-1$} or simply a \textit{great
     sphere.}
      \item For $A$ a closed subset of $S$, the set ${\rm Clos}(x,A)$ of
     points in $A$ closest to a point $x\in S$ is given by
     \[
    {\rm Clos}(x, A)=\{ a\in A \, |\, \forall b\in A, \ |a-x|\le |b-x| \} \, .
     \]
      \item $\mathcal{V}=\{V_1,V_2,\dots, V_s\}$ denotes a collection of
    codimension $1$ affine subspaces of $\mathbb{R}^N$ containing the
    center of $S$. For any $x$ in $S$, $\mathcal{V}_x$ will stand for the
    collection of $V_i$'s that contain $x$.
    \item For each $\alpha$ a subset of $\{1,2,\dots,s\}$ we define the subspace
    \[
    V_\alpha=\bigcap_{i \in \alpha} V_i
    \]
    and the great sphere in $S$
    \[
    S_\alpha=S \cap V_\alpha\,.
    \]
    By usual convention, $V_\emptyset =\mathbb{R}^N$ and $S_\emptyset =S.$
    \item We call the pair $\mathbb{S}=(S,\mathcal{V})$ a
    \textit{stratification} of $S$ induced by the subsets of $\mathcal{V}.$
    We call the lowest dimensional $S_\alpha$ containing a point $x$
    the \textit{stratum of} $x$ and denote it by $S_\alpha(x).$ It is
    well-defined, since the intersection of any such two strata would be a lower
    dimensional one. A stratification induces a partition of $S$ into
    \textit{cells} of dimensions $0,1,\dots,\dim S$, the \textit{cell
    containing} $x$, denoted $\sigma(x),$  being
    the connected component of $\{ y\in S_\alpha(x)|\ S_\alpha(y)=S_\alpha(x) \}$
    containing $x.$ Two stratifications $\mathbb{S}_1=(S_1,\mathcal{V}_1)$ and
    $\mathbb{S}_2=(S_2,\mathcal{V}_2)$ are called \textit{isometric}
    if there exists a bijection
    $j:\mathcal{V}_1 \rightarrow \mathcal{V}_2$
    and an orthogonal affine map $O: S_1 \rightarrow S_2$ such that
    \[
    O(S_1\cap V_i)=S_2\cap j(V_i)
    \]
    for every $V_i \in \mathcal{V}_1.$
    \item Let $V_\alpha^\perp \subset \mathbb{R}^N$ be the affine subspace
    orthogonal to $V_\alpha \subset \mathbb{R}^N$ passing through the
    center of $S$ and $V_\alpha^\perp(\epsilon) \equiv
    {\rm U}_\epsilon(c)\cap V_\alpha^\perp$.
    \item The \textit{$\epsilon$-tubular neighborhood of $S_\alpha$} is
    defined as
    \begin{equation}\label{E:tubular nhd}
    T_\epsilon(S_\alpha)=(V_\alpha \times V_\alpha^\perp(\epsilon))\cap S\,.
    \end{equation}
    When the $\dim S_\alpha = k$ we also call $ T_\epsilon(S_\alpha)$ a
    \textit{k-tube}.
    \item For every $z\in S_\alpha$ we denote the great sphere
    orthogonal to $S_\alpha$ at $z$ by $S_\alpha^\perp(z).$ Notice that
    \[
    S_\alpha^\perp(z)=(\mathbb{R}\cdot(z-c)\times V_\alpha^\perp)\cap S\,.
    \]
      \item For every $z\in S_\alpha$ it is easy to check that the set
    \[
    F_\alpha(z,\epsilon)=(\mathbb{R}^+\cdot (z-c)\times
    V_\alpha^\perp(\epsilon))\cap S\,.
    \]
    is a spherical cap. We call it an
    $(\alpha, \epsilon)-$\textit{spherical cap with top at} $z$ or simply a
    \textit{spherical cap}.
    \item  For every $z\in S_\alpha$  we define a neighborhood fibered by
    $(\alpha, \epsilon)-$spherical caps
    \begin{equation}\label{E:patch}
    {\rm P}_\alpha(z,\epsilon)= \bigcup_{y\in
    {\rm U}_\epsilon(z)\cap\sigma(z)}F_{\alpha}(y,\epsilon )\ ,
    \end{equation}
    called the {\em $(\alpha,\epsilon)$-patch centered at} $z$
    where $\alpha$ is given by the stratum containing $z$.
    \item The boundary $\partial_{rel} F_\alpha(z,\epsilon)$ of the spherical cap
    $F_\alpha(z,\epsilon)$ relative to the subspace $S_\alpha^\perp(z)$ is given by
    \[
    \partial_{rel} F_\alpha(z,\epsilon)
    =\left(\mathbb{R}^+\cdot (z-c)\times
    (\partial {\rm U}_\epsilon(c)\cap V_\alpha^\perp) \right)\cap S \,.
    \]
    The relative boundaries of spherical caps are Euclidean spheres though not great ones: rather
    they are analogous to meridians with the north pole at $z$.
    Note that
    \begin{equation}\label{E:dim}
    \dim \partial_{rel} F_\alpha(z,\epsilon)= n- \dim S_\alpha-1
    =\dim S - \dim V_\alpha\,.
    \end{equation}

    \item The stratification $\mathbb{S}$ induces the stratification
    $(\partial_{rel} F_\alpha(z,\epsilon),\mathcal{V}_z)$
    on the sphere $\partial_{rel} F_\alpha(z,\epsilon).$

\end{itemize}

\medskip
 We leave to the reader to check the following relations:
    \begin{equation}\label{E:tubes}
    T_\epsilon(S_\alpha)=
    \bigcup_{z\in S_{\alpha}}F_{\alpha}(z,\epsilon )\,,
    \end{equation}
 otherwise speaking, the $(\alpha,\epsilon)$-spherical caps foliate
 the $\epsilon$-tubular neighborhoods of $S_\alpha$;
    \begin{equation}
    F_\alpha(z,\epsilon)\subset {\rm U}_{(\sqrt{2})\epsilon}(z)\,;
    \end{equation}
    \begin{equation}\label{E:patchnhd}
       {\rm P}_\alpha(x,\epsilon)   \subset {\rm
       U}_{(1+\sqrt{2})\epsilon}(x)\,;
    \end{equation}
    \begin{equation}\label{E:spherepartition}
    S=\bigcup_{z\in S_{\alpha}}S^{\perp}_\alpha (z)\,;
    \end{equation}
and
    \[
    \partial_{rel} F_\alpha(z,\epsilon)=S(c+\frac{\sqrt{r^2-\epsilon^2}}{r}(z-c),\epsilon)\cap
    S_\alpha^\perp(z)\ .
    \]
We also leave as an exercise the following.
\begin{lemma}\label{L:capcharacterization}{\bf (Spherical Cap Charactization)}
A spherical cap is the intersection of a half-space and a sphere,
and conversely.
\end{lemma}
\subsection{Spherical caps, the Set $\boldsymbol{\Omega}$, and the Main Lemma}

\bigskip
\begin{lemma}\label{L:lemma2}{\bf(Spherical Cap Boundary Isometry)}
If $\mathcal{V}_x=\mathcal{V}_y$ and $x,y \in S_\alpha(x),$ then
$(\partial_{rel} F_\alpha(x,\epsilon),\mathcal{V}_x)$ and
$(\partial_{rel}F_\alpha(y,\epsilon),\mathcal{V}_y)$ are isometric
stratifications for any $\epsilon>0.$
\end{lemma}

\begin{proof}
Consider the plane of $x, y$ and $c.$ This plane is a subspace of
$V_\alpha.$ Then the affine orthogonal map $O,$ which rotates this
plane about $c$ rotating $y$ to $x$ and which is the identity on
the space orthogonal the plane, is the required isometry . $\Box$
\end{proof}
\bigskip
\begin{lemma}\label{L:lemma3} {\bf (Spherical Cap Intersection)}
 Let $S_{\alpha'}\subset S_{\alpha}.$
 Given $\mu < r,$ if  $x\in S_{\alpha},$ and $z \in S_{\alpha'},$ then
 either
\begin{align}\label{E:capintersection}
   S^\perp_\alpha(x)\cap F_{\alpha'}(z,\mu) &= \emptyset\,, \notag\\
\intertext{or} S^\perp_\alpha(x)\cap F_{\alpha'}(z,\mu)&=
   F_{\alpha}(x,\delta)\,,
\end{align}
where $\delta$ and $\mu$ are related by
\begin{equation}\label{E:cap2}
r^2-\mu^2 = (r^2-\delta^2)\left( \frac{(x-c)\cdot
(z-c)}{r^2}\right) ^2\ ,
\end{equation}
with $r$ standing for the radius of the sphere $S$ centered at
$c$.

Moreover, if $\frac{\delta}{\mu} \le \sqrt{2\sqrt{2}-2}\simeq
.911,$ and $\frac{\delta}{\mu} \le \tau$ then
 \begin{equation}\label{E:luneestimate}
      d(x,\partial_{rel} F_{\alpha'}(z,\mu))
        < \tau \delta \ .
 \end{equation}
\end{lemma}

\begin{proof} Observe that if $x\notin F_\alpha(x,\mu)$ then
$S^\perp_\alpha(x)\cap F_{\alpha'}(z,\mu)= \emptyset.$
\begin{figure}[ht]
 \centerline{\epsfbox{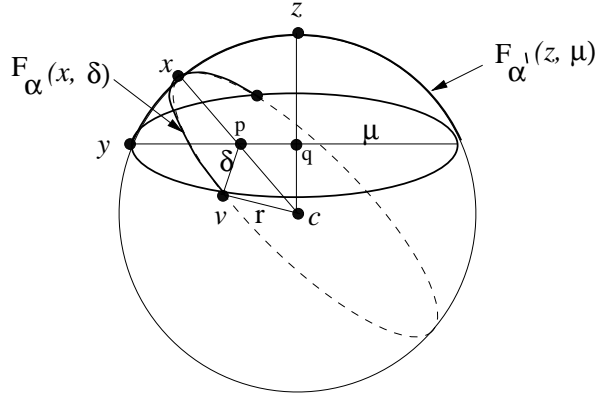}}
 \caption{Spherical cap intersecting a great 2-sphere }
 \label{F:cap}
\end{figure}
Now let $x\in F_\alpha(x,\mu)\setminus \partial_{rel}
F_{\alpha'}(z,\mu)$ and $x\ne z.$ then $S^\perp_\alpha(x)\cap
F_{\alpha'}(z,\mu)\ne \emptyset.$ Then $S^\perp_\alpha(x)\cap
F_{\alpha'}(z,\mu)\ne \emptyset.$ By Lemma
\ref{L:capcharacterization} $S^\perp_\alpha(x)\cap
F_{\alpha'}(z,\mu)$ is a shperical cap and we are justified in
denoted it by $F_\alpha(x,\delta).$ We must now determine
$\delta.$

Choose $v\in \partial F_{\alpha'}(z,\mu).$ Then $v\in \partial_
{rel} F_{\alpha}(x,\delta)\cap
\partial_{rel} F_{\alpha'}(z,\mu).$ Then the vectors
$v-c, \ x-c, \ z-c$ determine the three
dimensional $c$-centered sphere $S_3$ that contains the points
$v,x,$ and $z.$ We assume that $c$ is at the origin and the point
$z$ is at the north pole. The intersection of $F_{\alpha'}(z,\mu)$
with $S_3$ is a spherical cap whose boundary is a meridian which
is also $\partial_{rel} F_{\alpha'}(z,\mu)\cap S_3;$ this meridian
is determined by $\mu$ and the intersection of $\partial_{rel}
F_{\alpha}(x,\delta)$ with $S_3.$ This intersection is an arc of a
great circle whose endpoints lie on this meridian. Let $y$ be the
intersection of the great circle determined by $x$ and $z$ with
this meridian. The projections of $v$ on $x$ and $z$ are
$p=(\frac{v\cdot x}{r^2})x$ and $q=(\frac{v\cdot z}{r^2})z$
respectively. Note that
\begin{equation}\label{E:cos1}
\frac{v\cdot z}{r^2} = \sqrt{1-(\frac{\mu}{r})^2}\,,
\end{equation}
and that
\begin{equation}\label{E:cos2}
   \frac{v\cdot x}{r^2}=\sqrt{1-(\frac{\delta}{r})^2} \ .
\end{equation}
We can write
\begin{equation}\label{E:v}
     v=(\frac{v\cdot x}{r^2})x + u \,,
\end{equation}
which exhibits the contribution $u\in V^\perp_\alpha \subset
V^\perp_{\alpha'}$ to $v$ and yields
\begin{equation}\label{E:vzDot z}
v\cdot z= (\frac{v\cdot x}{r^2})x\cdot z \ .
\end{equation}
Equation \eqref{E:cap2} then follows by substituting
\eqref{E:cos1} and \eqref{E:cos2} in \eqref{E:vzDot z}.

\medskip
For the case $x=z$ it is apparent that \eqref{E:cap2} holds with
$\delta = \mu.$

\medskip
 Finally, for the case $x\in\partial_{rel}
F_{\alpha'}(z,\mu)$ it is easy to see that \eqref{E:cap2} holds
with $\delta = 0.$

\medskip
To obtain Equation \eqref{E:luneestimate} we first prove that $y$
is the closest vector to $x$ in $\partial_{rel}
F_{\alpha'}(z,\mu).$ Let $v$ play the role of any vector in
$\partial_{rel}F_{\alpha'}(z,\mu),$ not just one of those at a
distance $\delta$ from the axis of $S$ through $x,$ and observe
that $d(x,v)> d(x,y).$

\medskip
From the fact that angle xpy is acute--\textit{i.e.,}
\begin{equation}
  d(x,y)^2 \le d(y,p)^2 + d(p,x)^2\,,
\end{equation}
and the fact that the shortest distance to the circumference is
along a radius--\textit{i.e.,}
\begin{equation}
  d(p,x) \le  d(p,y)\,,
\end{equation}
we get
\begin{equation}
d(x,y)\le \sqrt{2}d(y,p)\,.
\end{equation}
Furthermore, from
\[
d(y,p)=d(y,q)-d(p,q)\,,
\]
and the fact that $\bigtriangleup vpq$ is a right triangle and we
have that
\begin{align}\label{E:delta}
 d(x,y) &\le \sqrt{2}(\mu - \sqrt{\mu^2 -\delta^2})\\
 &= \sqrt{2}\mu\left(1-\sqrt{1-\left(\frac{\delta}{\mu}\right)^2}\right) \,.\notag
\end{align}
A straight forward computation shows
$\frac{\delta}{\mu}\le\sqrt{2\sqrt{2}-2}$ is equivalent to
\begin{equation}
1-\sqrt{1-\left(\frac{\delta}{\mu}\right)^2}\le
\frac{1}{\sqrt{2}}\left(\frac{\delta}{\mu}\right)^2\,.
\end{equation}
Thus
\[
d(x,y)\le \mu\left(\frac{\delta}{\mu}\right)^2\le \tau \delta\,.
\]

\medskip
For the case $x=z$ it is apparent that \eqref{E:capintersection}
holds with $\delta = \mu.$

\medskip
 Finally, for the case $x\in\partial_{rel}
F_{\alpha'}(z,\mu)$ it is easy to see that \eqref{E:cap2} holds
with $\delta = 0.$

$\Box$
\end{proof}

\bigskip
The following lemma is an immediate corollary of the first part of
the preceding one where $z$ is replaced by $z'$, $x$ is replaced
by $z,$ $\mu$ is replaced by $\epsilon',$ and $\epsilon$ is taken
to be a number not necessarily equal to $\delta.$ In fact Case I
corresponds to $\epsilon>\delta,$ Case II to $\epsilon=\delta,$
and Case III to either $\epsilon<\delta$ or an empty intersection.
\begin{lemma}\label{L:capintersection} Let $S_{\alpha'}\subset
S_{\alpha}.$ Given $\ z \in S_{\alpha},$ and $z' \in S_{\alpha'},$
one of the following mutually exclusive relations hold:
\begin{align}
 \text{{\rm I}.  }&S^\perp_\alpha(z)\cap F_{\alpha'}(z',\epsilon')
                   \supsetneq F_{\alpha}(z,\epsilon); \notag\\
 \text{{\rm II}.  }&S^\perp_\alpha(z)\cap F_{\alpha'}(z',\epsilon')=
                    F_{\alpha}(z,\epsilon); \notag\\
 \text{{\rm III}. }&S^\perp_\alpha(z)\cap F_{\alpha'}(z',\epsilon')
 \subsetneq F_{\alpha}(z,\epsilon)\notag \, .
\end{align}

\end{lemma}

\bigskip
\begin{numberedremark}\label{R:capradius} Observe that by the
Pythagorean Theorem $d(z,\partial_{rel} F_{\alpha'}(z,\mu)))>\mu.$
\end{numberedremark}

\begin{figure}[ht]
 \centerline{\epsfbox{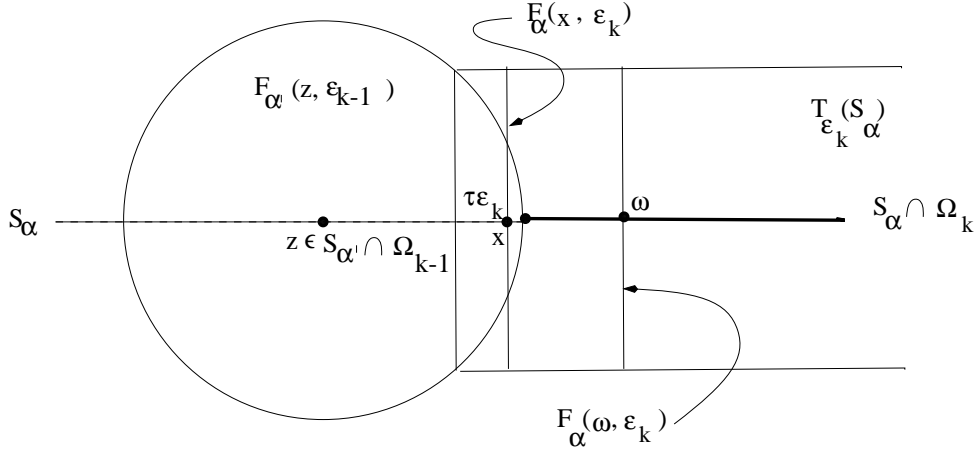}}
 \caption{Schematic aid }
 \label{F:aid}
\end{figure}

The authors found schematic representations of Figure \ref{F:cap}
and Figure \ref{F:aid} useful mental aids for arguments to follow.

\bigskip
\begin{definition}\label{D:Omega} \textbf{(The set $\boldsymbol{\Omega})$}
With respect to an $n$-tuple
$\boldsymbol{\epsilon}=(\epsilon_0,\epsilon_1,\dots,\epsilon_{n-1})$
of positive numbers set
\[
 \mathbf{T}_j=\mathbf{T}_j(\epsilon_j)=\bigcup_{\dim S_\alpha=j} T_{\epsilon_j}(S_\alpha)\,,
\]
for $j=0,1,\dots,n-1.$  Next set
\[
 \Omega_0=(\bigcup_{\dim S_\alpha=0} \ S_\alpha)
\]
and
\[
  \Omega_k=\Omega_k(\epsilon_{0}, \dots, \epsilon_{k-1})=(\bigcup_{\dim S_\alpha=k} \
  S_\alpha)\setminus \bigcup_{j<k} \mathbf{T}_j(\epsilon_j)\,.
\]
for $k=1,\dots,n.$ Note that sets in this hierarchy may be empty
because there are not enough affine subspaces in $\mathcal{V}$ or
the $\epsilon_j$ are too big. Finally, set
\[
  \Omega=\Omega(\boldsymbol{\epsilon})=\bigcup_{k\le n} \Omega_k\,.
\]
\end{definition}

\begin{definition}\label{D:epsilon}
 An n-tuple
$\boldsymbol{\epsilon}=(\epsilon_0,\epsilon_1,\dots,\epsilon_{\dim
S-1})$ of positive numbers is said to \textit{decrease
sufficiently rapidly and prevent improper tube interference} with
respect to a stratification $(S, \mathcal{V})$ if

\noindent(1) for $k>0,\ \epsilon_k / \epsilon_{k-1}\le \tau,
\text{ with } \epsilon_0 \le \tau r,$ where $\tau
=\frac{\sqrt{2-2\sqrt{\frac{2}{3}}}}{\ 8}\,,$

\noindent(2) for $0\le k\le\dim S -1$ and $\omega\in \Omega_k,$ if
$V_\alpha \subset V_i \notin \mathcal{V}_\omega , $ then ${\rm
U}_{8\epsilon_k}(\omega)\cap T_{\epsilon_k}(S_\alpha)=\emptyset
\,.$

\end{definition}

Item (1) will ensure that the $\epsilon_k$'s decrease fast enough.
The value of $\tau$ is chosen so that
\begin{equation}\label{E:tau}
{\rm U}_{8\epsilon_k}(z) \subset {\rm
U}_{\left(\sqrt{2-2\sqrt{\frac{2}{3}}}\,\right)\epsilon_{k-1}}(z)\,.
\end{equation}
This relation implies that a $8\epsilon_k$-neighborhood of any
point in a sphere of radius $\epsilon_{k-1}$ is contained in an
octant of that sphere centered around the given point.  In
subsection 5.4.2, claim 2, of the proof of the Main Lemma we shall
need that the hypotenuse of a geodesic spherical right triangle is
longer than the other sides. This will be the case if the triangle
lies in a neighborhood contained in an octant. In addition this
value of $\tau$ is smaller than the quantity $\sqrt{2\sqrt{2} -2}$
which appears in Lemma \ref{L:lemma3}. It is also small enough so
that each time a geodesic belongs to the discussion (\textit{e.g.}
see subsection 5.4.3, claim 3 in the proof of the Main Lemma) it
is unique.

Item (2), which might impose that ratios of successive
$\epsilon_k$ are much smaller than $\tau,$ guarantees that
spherical caps fibering tubes in general do not interfere in
unexpected ways with spherical caps centered on points of $\Omega$
and leaves enough room to simplify some arguments. We call
attention to the reader that the subscript $k$ in
$T_{\epsilon_k}(S_\alpha)$ is not necessarily the same dimension
as that of $S_\alpha$ as it would be in the definition of
$\mathbf{T}_k$ in Definition \ref{D:Omega}.

\begin{lemma}\label{L:maximal}
With respect to the stratification $(S,\mathcal{V})$ let
$\boldsymbol{\epsilon}=(\epsilon_0,\epsilon_1,\dots,\epsilon_{n-1})$
decrease sufficiently rapidly and prevent improper tube
interference. Then for each $x\in \bigcup_{j<n} \mathbf{T}_j$
there exist a largest integer $ 0\le k(x)\le n-1,$ a great sphere
$S_\alpha$ with $\dim{S_\alpha}=k(x),$ and a point $z\in S_\alpha$
such that the $(\alpha,\epsilon_{k(x)})$-spherical cap
$F_\alpha(z,\epsilon_{k(x)})$ contains $x$  and satisfies
\[
\partial_{rel} F_\alpha
(z,\epsilon_{k(x)})\cap(\mathbf{T}_{-1}\cup \mathbf{T}_1 \cup
\dots \mathbf{T}_{k(x)-1})=\emptyset\,,
\]
 where
$\mathbf{T}_{-1}=\emptyset.$

Moreover, this $(\alpha,\epsilon_{k(x)})-$spherical cap with top
at $z$ and containing $x$ is unique. It will be denoted by
$F_{\alpha(x)}(z(x),\epsilon_{k(x)})$ to show all the dependencies
on $x,$ however, for the sake of brevity we shall also denote it
simply by $F_{\alpha}(z,\epsilon_{k(x)}).$ In the abbreviated form
it is to be understood that the dependence of $\alpha$ and $z$ on
$x$ is adequately indicated by the subscript $k(x)$ on $\epsilon.$
Here the unique great sphere $S_\alpha$ containing $z$ determined
by $x$ we denote by $S_{\alpha(x)}.$
\end{lemma}

\medskip
\begin{proof}
Take the largest $k_1 < n$ such that $x\in \mathbf{T}_{k_1}.$ By
\eqref{E:tubes} $x$ belongs to some spherical cap
$F_{\alpha_1}(z_1,\epsilon_{k_1})$ where $\dim S_{\alpha_1} =
k_1.$ If $\partial_{rel} F_{\alpha_1}(z_1,\epsilon_{k_1})\cap
\mathbf{T}_j = \emptyset$ for all $j<k_1$ we are done and the
spherical cap $F_{\alpha_1}(z_1,\epsilon_{k_1})$ is unique.

The argument for the uniqueness is the following.  There are two
cases:
\begin{enumerate}
\item $z_1\in \Omega_{k_1}\cap T_{\epsilon_{k_1}}(S_{\alpha_1})\,;$
\item  $z_1\in T_{\epsilon_{j}}(S_{\alpha'}),j<k_1, \ \dim
S_{\alpha'}=j, \text{ and } S_{\alpha'}\subset S_{\alpha_1}.$
\end{enumerate}
In case (1) $F_{\alpha_1}(z_1,\epsilon_{k_1})\subset {\rm
U}_{8\epsilon_{k_1}}(z_1).$ The uniqueness follows from the fact
that there are no improper tube interferences.

In case (2) we have $z_1\in T_{\epsilon_{j}}(S_{\alpha'})$ with
$j<k_1,$ say $z_1\in F_{\alpha'}(z,\epsilon_j).$ Lemma
\ref{L:lemma3} implies that $d(z_1,\omega)\le \tau\epsilon_{k_1}$
for some $\omega\in \Omega_{k_1}$ (in that lemma the role of $x$
is played by $z_1,$ $\delta$ by $\epsilon_{k_1},$ $\mu$ by
$\epsilon_{j},$ and $z$ by itself. Here
$F_{\alpha_1}(z_1,\epsilon_{k_1})\subset {\rm
U}_{8\epsilon_{k_1}}(\omega).$

Hence also in case (2), the uniqueness follows by the way
$\epsilon_{k_1}$ is chosen so that
$F_{\alpha_1}(z_1,\epsilon_{k_1})$ does not intersect any other
$(\epsilon_{k_1},\alpha'')$-spherical cap with $\dim S_{\alpha''}
= k_1.$

If $\partial_{rel} F_{\alpha_1}(z_1,\epsilon_{k_1})\cap
\mathbf{T}_j \neq \emptyset$ for some $k_2<k_1,$  take the largest
such $k_2.$ By Lemma \ref{L:capintersection} and definition of
$\epsilon_{k_1}$ we get that $x\in
F_{\alpha_1}(z_1,\epsilon_{k_1})\subset \mathbf{T}_{k_2}$ so that
$x$ belongs to a spherical cap $F_{\alpha_2}(z_2,\epsilon_{k_2})$
where $\dim S_{\alpha_2} =k_2.$ Either we are done and this cap is
unique or we must repeat the procedure finding $k_3$ and so forth.
This procedure terminates at some $k\ge 0$ which defines $k(x)$
and a unique spherical cap $F_\alpha (z,\epsilon_{k(x)})$
satisfying the conditions of this lemma. $\Box$
\end{proof}

\bigskip
\begin{mainlemma}
Let $(S,\mathcal{V})$ be a stratification with $\dim S=n$ and
$\Omega(\boldsymbol{\epsilon})=\Omega(\epsilon_0,\dots,\epsilon_{n-1})$
as in Definition \ref{D:Omega}. If
$\boldsymbol{\epsilon}=(\epsilon_0,\epsilon_1,\dots,\epsilon_{n-1})$
decreases sufficiently rapidly and prevents improper tube
interference, then
\[
{\rm clos}(x,\Omega)\subset \overline{\sigma(x)}\, .
\]
\end{mainlemma}

Before giving a proof of this lemma, we need some more notation.

\subsection{Induced Stratification and Related Constructs} Choose some $\alpha$ and set
$\dim V_\alpha = k+1$, hence $\dim S_\alpha =k \, .$ Pick any
$x\in S_\alpha$, and set $\widehat{S}=\partial_{rel}
F_\alpha(x,\epsilon_k)$ and $\widehat{\mathcal{V}}=\mathcal{V}_x=
\{V_i|\ i\in \alpha\}\,.$ Recall (see \eqref{E:dim}) that the
sphere $\widehat{S}$ has dimension $n-k-1.$ We denote the
stratification induced by $\widehat{\mathcal{V}}$ on
$\partial_{rel} F_\alpha(x,\epsilon_k)$ by
$\widehat{\mathbb{S}}=(\widehat{S},\widehat{\mathcal{V}}).$
\medskip
Let
\begin{align}\label{E:epsilon hat}
  \widehat{\boldsymbol{\epsilon}}&=(\widehat{\epsilon}_0,\dots,\widehat{\epsilon}_{n-k-2})\,, \\
         &=(\epsilon_{k+1},\dots,\epsilon_{n+k-1})\,. \notag
\end{align}
For $\breve {\alpha} \subset \alpha$ we denote the great sphere
$\widehat{S}\cap V_{\breve {\alpha}}$ by $\widehat{S}_{\breve
{\alpha}}$ and the $\epsilon$-tubular neighborhood $(V_{\breve
{\alpha}}\times V_{\breve {\alpha}}^\perp(\epsilon))\cap
\widehat{S}$ of $\widehat{S}_{\breve {\alpha}}$ in $\widehat{S}$
by $\widehat{T}_{\epsilon}(\widehat{S}_{\breve {\alpha}}) \, .$
Notice that
\begin{equation}\label{E:S_hat}
 \widehat{S}_{\breve {\alpha}}= S_{\breve {\alpha}}\cap \widehat{S}\,,
\end{equation}
and that by \eqref{E:tubular nhd}
\begin{equation}\label{E:tube hat}
 \widehat{T}_\epsilon(\widehat{S}_{\breve {\alpha}})=
 T_\epsilon(S_{\breve {\alpha}})\cap \widehat{S}\,.
\end{equation}
Furthermore, one can check that the cells $\widehat{\sigma}(x)$ of
the induced stratification are obtained by intersecting the
original cells with $\widehat{S}$: \textit{i.e.},
\begin{equation}\label{E:cell hat}
\widehat{\sigma}(x)=\sigma(x)\cap\widehat{S}\,.
\end{equation}
Next set
\begin{equation}\label{E:T hat}
 \widehat{\mathbf{T}}_j(\widehat{\epsilon}_{j})=
 \bigcup_{\substack{\dim \widehat{S}_{\breve {\alpha}}=j\\
 \breve {\alpha}\subset \alpha}}
 \widehat{T}_{\widehat{\epsilon}_j}(\widehat{S}_{\breve {\alpha}})\,,
\end{equation}
for $j=0,1,\dots,n-k-2,$ where by \eqref{E:epsilon hat} and
\eqref{E:tube hat}:
\[
\widehat{T}_{\widehat{\epsilon}_j}(\widehat{S}_{\breve {\alpha}})=
{T}_{\epsilon_{j+k+1}}(S_{\breve {\alpha}})\cap \widehat{S}\, .
\]
From the obvious relation
\[
\dim V_{\breve {\alpha}}=\dim V_\alpha + \dim( V_{\breve
{\alpha}}\cap V_{\alpha}^\perp)
\]
and the definitions of $\widehat{S}$ and of $j$ (as $\dim
\widehat{S}_{\breve {\alpha}}$), we get $\dim V_{\breve
{\alpha}}=j+k+2 \, ,$ thus
\begin{align}\label{E:bftubehat}
\widehat{\mathbf{T}}_j(\widehat{\epsilon}_{j})&=
 \bigcup_{\substack{\dim S_{\breve {\alpha}}=j+k+1 \\
 \breve {\alpha}\subset \alpha}}
 T_{\epsilon_{j+k+1}}(S_{\breve {\alpha}})\cap \widehat{S} \notag\\
&=\mathbf{T}_{j+k+1}(\epsilon_{j+k+1})\cap \widehat{S}\, .
\end{align}

Next set
\[
  \widehat{\Omega}_0=\bigcup_{\substack{\dim \widehat{S}_{\breve {\alpha}}=0\\
  \breve {\alpha}\subset \alpha}} \ \widehat{S}_{\breve {\alpha}}
\]
and
\[
  \widehat{\Omega}_m=
  \widehat{\Omega}_m(\widehat{\epsilon}_{0}, \dots, \widehat{\epsilon}_{m-1})=
  (\bigcup_{\substack{\dim \widehat{S}_{\breve {\alpha}}=m \\
  \breve {\alpha}\subset \alpha}} \
  \widehat{S}_{\breve{\alpha}})\setminus \bigcup_{j<m}
  \widehat{\mathbf{T}}_j(\widehat{\epsilon}_{j})\,.
\]
for $m=1,\dots,n-k-1.$ Finally, set
\[
 \widehat{\Omega}=\widehat{\Omega}(\widehat{\boldsymbol{\epsilon}})=
 \bigcup_{m\le n-k-1} \widehat{\Omega}_m\, .
\]
From \eqref{E:S_hat} and \eqref{E:bftubehat} we obtain the key
relation:
\begin{equation}\label{E:Omega}
\widehat{\Omega}=\Omega \cap \widehat{S}\,.
\end{equation}

\bigskip
\subsection{Proof of Main Lemma}

We shall prove this by induction on the dimension of S.

\begin{hypothesis}
Let $(S,\mathcal{V})$ be a stratification with $\dim{S} <n.$ If
$\boldsymbol{\epsilon}=(\epsilon_0,\epsilon_1,\dots,
\epsilon_{\dim{S-1}})$ decreases sufficiently rapidly and prevents
improper tube interference, then for $x\in S$
\[
 {\rm clos}(x,\Omega)\subset
 \overline{\sigma(x)} \, .
 \]

\end{hypothesis}

After some preliminary discussion, the proof will consist in
proving four claims which localize ${\rm clos}(x,\Omega)$ with
greater and greater precision. The fourth and final claim finishes
the induction and hence the proof of the  Main Lemma 5.6.

\medskip
 Repeating here the contents of Definition \ref{D:epsilon} for the
 convenience of the reader, the premise of the inductive hypothesis states that
 \begin{enumerate}
   \item for $m>0, \ \epsilon_m/ \epsilon_{m-1}\le \tau,$
  with $\epsilon_0 \le \tau r, $ where $\tau = \frac{\sqrt{2-2\sqrt{\frac{2}{3}}}}{\ 8}\,;$
  \item  for $0\le m \le \dim S-1$ and $\omega\in
      \Omega_m,$ if $V_{\alpha} \subset
   V_i \notin \mathcal{V}_{\omega}\, ,$ then
  ${\rm U}_{8\epsilon_m}(\omega)\cap
  T_{\epsilon_m}(S_{\alpha})=\emptyset\,.$

\end{enumerate}
For a stratification $(\widehat{S},\widehat{\mathcal{V}})$ induced
by the stratification $(S, \mathcal{V})$ we get that
$\widehat{\epsilon}_m$ inherits property (1) from
$\epsilon_{m+k+1}.$ Also ${\rm
U}_{8\widehat{\epsilon}_m}(\widehat{\omega})$ inherits property
(2) from ${\rm U}_{8\epsilon_{m+k+1}}(\widehat{\omega}):$ \textit{
i.e.,} ${\rm U}_{8\widehat{\epsilon}_m}(\widehat{\omega})\cap
\widehat{T}_{\widehat{\epsilon}_m}(\widehat{S}_{\widehat{\alpha}})=\emptyset
\,, $ since ${\rm U}_{8\epsilon_{m+k+1}}(\widehat{\omega})\cap
T_{\epsilon_{m+k+1}}(S_{\widehat{\alpha}})=\emptyset$ whenever
$V_{\widehat{\alpha}} \subset V_i\notin
\mathcal{V}_{\widehat{\omega}}\,.$

Under the assumption that the premise of the inductive hypothesis
holds for the stratification $(S,\mathcal{V})$ where $\dim S<n,$
we get from the above inheritance properties that the premise
holds for all induced stratifications
$(\widehat{S},\widehat{\mathcal{V}})$ with $\dim \widehat{S}<n-1.$
So the conclusion of the inductive hypothesis then holds for these
lower dimensions. This fact will come into play in the final
claim.

\bigskip
To start the induction at $n=1$, observe that when $\dim{S}=0,\ S$
consists of two points and the induction hypothesis (1)
immediately follows.

\bigskip
Next assume that the inductive hypothesis(n) holds. To prove the
induction hypothesis (n+1) let $S$ be a sphere $\dim{S}=n$ and
$\boldsymbol{\epsilon}=(\epsilon_0,\epsilon_1,\dots,
\epsilon_{\dim{S-1}})$ decrease sufficiently rapidly and prevent
improper tube interference.

\medskip
Take $x\in S$. If $x\in \Omega$ there is nothing to prove. Suppose
$x\notin \Omega$. Then $x\in \bigcup_{j<n} \mathbf{T}_j,$ and we
get a unique spherical cap $F_{\alpha}(z,\epsilon_{k(x)})$ which
satisfies Lemma \ref{L:maximal}. Recall Equation \eqref{E:patch}:
\textit{i.e.},
\[
{\rm P}_{\alpha}(z,3\epsilon_{k(x)})=\bigcup_{y\in {\rm
U}_{3\epsilon_{k(x)}}(z)\cap \sigma(z)}
F_{\alpha}(y,3\epsilon_{k(x)})\,.
\]
(In this expression, as in the convention for spherical caps, the
dependence of $\alpha$ and $z$ on $x$ is indicated by the
subscript $k(x)$ on $\epsilon.$) Either $z=\omega\in
\Omega_{k(x)},$ or $z\notin \Omega_{k(x)}$ and there exists
$\omega \in \Omega_{k(x)}$ such that $d(z,\omega)< \tau
\epsilon_{k(x)}$ by using Lemma \ref{L:lemma3}. In either case by
\eqref{E:patchnhd} we have
\begin{equation}\label{E:patchnhd2}
{\rm P}_{\alpha}(z,3\epsilon_{k(x)})\subset {\rm
U}_{8\epsilon_{k(x)}}(\omega) \,.
\end{equation}

\subsubsection {Claim 1:} ${\rm clos}(x,\Omega)\subset
{\rm P}_{\alpha}(z,3\epsilon_{k(x)}).$

\noindent\textit{Proof:} Due to the definition of $k(x)$ and
property (2) of Definition \ref{D:epsilon}, the sphere
$\partial_{rel} F_{\alpha}(z,\epsilon_{k(x)})$ has a non empty
intersection with $\Omega$: so the distance from $x\in
F_{\alpha}(z,\epsilon_{k(x)})$ to $\Omega$ is less than or equal
to the ${\rm diameter}(F_{\alpha}(z,\epsilon_{k(x)})) =
2\epsilon_{k(x)}.$ \textit{Claim 1} follows because
$d(x,S\setminus {\rm
P}_{\alpha}(z,3\epsilon_{k(x)}))>2\epsilon_{k(x)}.$

\subsubsection {Claim 2:} ${\rm clos}(x,\Omega)\subset
F_{\alpha}(z,3\epsilon_{k(x)})\cup ({\rm
U}_{3\epsilon_{k(x)}}(z)\cap \sigma(z))$.

\noindent\textit{Proof:} Let $\omega \in {\rm clos}(x,\Omega).$ We
know by \textit{Claim 1} that $ \omega \in {\rm
P}_{\alpha}(z,3\epsilon_{k(x)})):$ so we can assume that $\omega
\in {\rm P}_{\alpha}(z,3\epsilon_{k(x)})\setminus {\rm
U}_{3\epsilon_{k(x)}}(z)\cap\sigma(z)$. We just have to prove that
$\omega \in F_{\alpha}(z,3\epsilon_{k(x)}). $  From \textit{Claim
1} there exists $y\in {\rm U}_{3\epsilon_{k(x)}}(z)\cap \sigma(z)$
such that $\omega\in F_{\alpha}(y,3\epsilon_{k(x)})\setminus {\rm
U}_{3\epsilon_{k(x)}}(z)\cap \sigma(z).$ Project the point
$\omega$ to $F_{\alpha}(z,3\epsilon_{k(x)})$ by using the isometry
given in Lemma \ref{L:lemma2}, say
\[
\omega'=O^{-1}(\omega)\in F_{\alpha}(z,3\epsilon_{k(x)})\setminus
F_{\alpha}(z,\epsilon_{k(x)})\,.
\]
Now there are two possibilities.
\[
 \text{Either }\omega'\in \Omega, \text{ or } \omega'\in \mathbf{T}_j, j<{k(x)}\, .
\]
This last possibility does not occur because it would imply
\[
\partial_{rel} F_{\alpha}(z,\epsilon_{k(x)})\subset \mathbf{T}_j\,,
\]
by using Lemma \ref{L:lemma3}, contradicting the definition of
${k(x)}$ in Lemma \ref{L:maximal}. Hence, $\omega'\in \Omega$.

Assume $\omega'\ne \omega.$ Then by the property of geodesic right
triangles resulting from Definition \ref{D:epsilon}(1) we get
$|\omega-x|>|\omega'-x|$ which contradicts the fact that
$\omega\in {\rm clos}(x,\Omega).$ Hence $\omega=\omega'.$

\subsubsection {Claim 3:} ${\rm clos}(x,\Omega)\subset
\partial_{rel}F_{\alpha}(z,\epsilon_{k(x)})\cup ({\rm
U}_{3\epsilon_{k(x)}}(z)\cap \sigma(z)).$

\noindent\textit{Proof:} Let $\omega \in {\rm clos}(x,\Omega).$ We
know by \textit{Claim 2} that $ \omega \in
F_{\alpha}(z,3\epsilon_{k(x})\cup ({\rm
U}_{3\epsilon_{k(x)}}(z)\cap \sigma(z)):$ so we can assume
$\omega\in  F_{\alpha}(z,3\epsilon_{k(x)})\setminus \{z\}$. Along
the great circle on $S$ from $\omega$ through $x$ we have a point
$\omega'\in\partial_{rel} F_{\alpha}(z,\epsilon_{k(x)})$ which is
well defined by taking the shortest geodesic from $\omega$ to
$\omega'.$ Again their are two possibilities.
\[
\text{Either }\omega'\in \Omega \text{ or } \omega'\in
\mathbf{T}_j, \text{ for some } j>{k(x)}\,.
\]
The last possibility does not occur for the following reason.
Using Lemma \ref{L:maximal} we can take $j=k(\omega')>k(x).$  This
means that $\omega '\in F_{\alpha(\omega
')}(z(\omega'),\epsilon_{k(\omega')})$ for some $z(\omega') \in
S_{\alpha(\omega')}\supset S_{\alpha}$ with $\dim
S_{\alpha(\omega')}=k(\omega')$ and $\partial_{rel}F(z(\omega'),
\epsilon_{k(\omega')})\cap (\mathbf{T}_0\cup \dots \cup
\mathbf{T}_{k(\omega')-1})=\emptyset.$ Also we would have that
$\omega\in clos(\omega',\Omega),$ for otherwise $\omega\notin
clos(x,\Omega)$.
\begin{figure}[ht]
 \centerline{\epsfbox{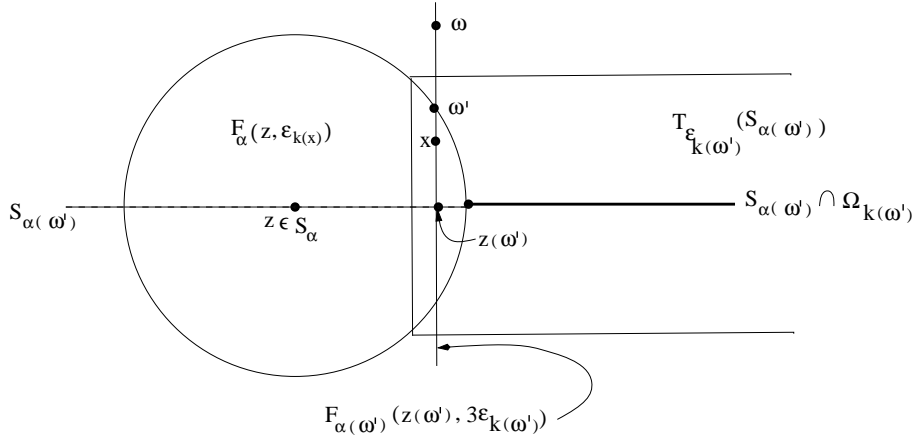}}
 \caption{Schematic Aid to \textit{Claim 3}.
 Here the disk corresponds to a transversal view of the horizontal strip in Figure \ref{F:aid}.}
 \label{F:aid2}
\end{figure}
Using \textit{claim 2}, with $\omega'$ playing the role of $x$ in
that claim, we would get $\omega\in
F_{\alpha(\omega')}(z(\omega'),3\epsilon_{k(\omega')}).$ Since
$\omega, \omega'\in S_{\alpha(\omega')}^\perp(z(\omega')),$ we
would also get $x\in S_{\alpha(\omega')}^\perp(z(\omega')).$
Applying Lemma \ref{L:capintersection} Case III, we would get
$x\in F_{\alpha(\omega')}(z(\omega'),\epsilon_{k(\omega')})$ (See
Figure \ref{F:aid2}) which would imply that $k(x)\ge k(\omega'),$
a contradiction completing the proof of \textit{claim 3}.

\subsubsection {Claim 4:} ${\rm clos}(x,\Omega)\subset
\overline{\sigma(x)}$.

\noindent\textit{Proof:} From the previous claim
\[{\rm clos}(x,\Omega)\subset \partial_{rel}
F_\alpha(z,\epsilon_{k(x)})\cup {\rm U}_{3\epsilon_{k(x)}}(z)\cap
\sigma(z)\, .
\]
If $x=z,$ then, by virtue of \eqref{E:luneestimate} and the fact
that the distance from $x$ to the boundary of the
$\epsilon_{k(x)}$-tubular neighborhood of $z$ is greater than
$\epsilon_{k(x)}$ (see Remark \ref{R:capradius}) we get,
\[
{\rm clos}(x,\Omega)\subset {\rm U}_{3\epsilon_{k(x)}}(z)\cap
\overline{\sigma(x)}\,.
\]
So assume $x\ne z$ and
\begin{equation}\label{E:assumption}
{\rm clos}(x,\Omega)\cap  \partial_{rel}
F_\alpha(z,\epsilon_{k(x)})\setminus {\rm
U}_{3\epsilon_{k(x)}}(z)\cap \sigma(z)\ne \emptyset \,.
\end{equation}

 Project $x\in F_\alpha(z,\epsilon_{k(x)})$ along the great circle on $S$ from $z$
through $x$ to $\widehat{x}\in \partial_{rel}
F_\alpha(z,\epsilon_{k(x)}).$ Since $x$ and $\widehat{x}$ lie on
an arc of a great circle which doesn't cross cell boundaries,
$\widehat{x}\in \overline{\sigma(x)}.$ Hence
\[
\overline{\sigma(x)}=\overline{\sigma(\widehat{x})}\,.
\]
If $\widehat{x}\in \Omega,$ then ${\rm clos}(x,\Omega)=\{\hat{x}
\}, $ and we are finished. If $\widehat{x}\notin \Omega,$ we use
the inductive hypothesis to get
\[
 {\rm clos}(\widehat{x},\widehat{\Omega})\subset \overline{\widehat{\sigma}(\widehat{x})}\, .
\]
We can assert for $\widehat{x}\in
\partial_{rel} F_\alpha(z,\epsilon_{k(x)})$ that
\begin{equation}\label{E:cellhat}
\overline{\widehat{\sigma}(\widehat{x})}=
\overline{\sigma(\widehat{x})} \cap
 \partial_{rel} F_\alpha(z,\epsilon_{k(x)}) \subset \overline{\sigma(\widehat{x})}
\end{equation}
despite the fact that it does not hold for spherical cap
boundaries in general. The reason for this is that affine
subspaces in $\mathcal{V}$ which intersect strata of
$\widehat{S}=\partial_{rel} F_\alpha(z,\epsilon_{k(x)})$ to form
cells in $\widehat{S}$ are precisely those in $\mathcal{V}$ which
contain $V_\alpha.$ From \eqref{E:Omega} and \eqref{E:assumption}
\[
 {\rm clos}(x,\Omega) ={\rm clos}(x,\widehat{\Omega})\, .
\]
However, we have
\[
 {\rm clos}(x,\widehat{\Omega}) ={\rm clos}(\widehat{x},\widehat{\Omega})\,,
\]
the reason for which is the following.

\begin{figure}[ht]
 \centerline{\epsfbox{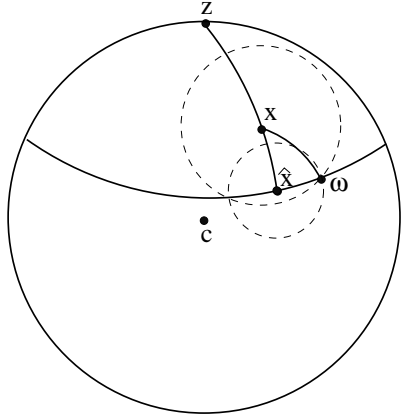}}
 \caption{Reduction to a great 2-sphere in three dimensions }
 \label{F:final}
\end{figure}

Observe that
\[
\widehat{S}\cap S(x,d(x,\Omega))= \widehat{S}\cap
S(\widehat{x},d(\widehat{x},\widehat{\Omega}))
\]
as illustrated in Figure \ref{F:final} where we have drawn the
great 2-sphere in three dimensions determined by the center $c$
and the points $z,x,\omega \in {\rm clos}(x,\widehat{\Omega})\, .$
$\widehat{S}$ intersects this sphere in the meridian containing
$\widehat{x}$ and $\omega.$ Thus
\begin{align}
  {\rm clos}(x,\widehat{\Omega})&=(\widehat{S}\cap S(x,d(x,\Omega)))\cap \widehat{\Omega}\notag \\
       &= (\widehat{S}\cap S(\widehat{x},d(\widehat{x},\widehat{\Omega})))\cap \widehat{\Omega}  \notag \\
       &={\rm clos}(\widehat{x},\widehat{\Omega}) \, . \notag
\end{align}

This finishes the proof \textit{claim 4 }and also the induction
step in the proof of the Main Lemma 5.6.
 $\Box$

\subsection{Conclusion of proof of Proposition \ref{proposition3}}
Let S=S(1,0) be the $N-1$ dimensional unit sphere centered at the
origin and $\mathcal{V}$ be the set of co-dimensional 1 subspaces
of $\mathbb{R}^N$ normal to the edges of the polytope $P.$

\subsubsection{Proof of Proposition \ref{proposition3}(3):
the closest points $\omega \in \Omega$ to $y \in \bar{\sigma}
\subset C_1(v_i)$ are in the same cell closure $\bar{\sigma}$.}
Each $C(v_i)$ is a polyhedral cone. The half-spaces, of which it
is the intersection, are determined by co-dimension 1 subspaces
orthogonal to the edges connecting $v_i$ to neighboring vertices.
Let $\mathcal{V}_P$ be the collection of co-dimension 1 subspaces
orthogonal to all the edges emanating from all the vertices of
$P$. Let the collection $\mathcal{V}$ of co-dimension 1 subspaces
of \S 5 contain $\mathcal{V}_P$. The Main Lemma holds for either
collection of co-dimension 1 subspaces because the partition of
$S$ into cells $\sigma$ resulting from $\mathcal{V}$ merely
refines the one resulting from $\mathcal{V}_P$. That is to say
that the cell containing ${\rm clos}(x,\Omega)$ given by
$\mathcal{V}$ is a subset of the cell containing ${\rm
clos}(x,\Omega)$ given by $\mathcal{V}_P.$ This fact is the basis
of the proof of Theorem \ref{T:byproduct2}.

All that is left of show is that there exists
$\boldsymbol{\epsilon}=(\epsilon_0,\epsilon_1,\dots,\epsilon_{n-1})$
which satisfies the premise of the Main Lemma 5.6: namely,
$\boldsymbol{\epsilon}$ decreases sufficiently rapidly and
prevents improper tube interference. Upon establishing this, the
conclusion of the Main Lemma 5.6 will then hold.

Let $x\in S.$ Because $\mathcal{V}$ is a finite collection, the
planes in $\mathcal{V}$ which do not pass through $x$ will have a
certain minimal distance to $x$. Hence, a small enough
neighborhood of $x$ avoids all planes that do not contain $x;$ and
from the definition of $\sigma(x),$ we then get that for every
$x\in S$ there is $\epsilon(x)>0$ such that the following holds:
\begin{enumerate}
    \item for all $y \in {\rm U}_{\epsilon (x)}(x)\cap \sigma(x)$,
    $\mathcal{V}_y=\mathcal{V}_x$\,;
    \item if  $\ V_\alpha \subset V_i \notin \mathcal{V}_x,$ then
    ${\rm U}_{8\epsilon(x)}(x)\cap T_{\epsilon (x)}(S_\alpha)= \emptyset \, .$

\end{enumerate}

By virtue of compactness of the various components of
$\Omega_k(\epsilon_{0},\dots,\epsilon_{k-1}),$ we can then
inductively choose a uniform $\epsilon_k$ for $\Omega_k$ so that:
\begin{enumerate}
 \item for $k>0,\ \epsilon_k / \epsilon_{k-1}\le \tau, \text{ with }
 \epsilon_0 \le \tau r,$ where $\tau = \frac{\sqrt{2-2\sqrt{\frac{2}{3}}}}{8}\,;$
  \item for every $\omega\in \Omega_k$ if $V_\alpha \subset V_i
  \notin
  \mathcal{V}_\omega , $ then ${\rm U}_{8\epsilon_k}(\omega)\cap
  T_{\epsilon_k}(S_\alpha)=\emptyset \,,$
\end{enumerate}
as needed.

\subsubsection{Proof of Proposition \ref{proposition3}(2):
for each cell $\sigma \subset C_1(v_i)$ there is a
$\epsilon_\sigma>0$ such that $\Omega\cap {\rm
U}_{\varepsilon_\sigma}(\bar{\sigma}) \subset \bar{\sigma}$ }
Since $d(z,\partial_{rel} F_\alpha(z,\epsilon_k))>\epsilon_k,$ we
can choose $\epsilon_\sigma=\epsilon_k,$ where $k$ is the
dimension of $\sigma.$

\subsubsection{Proof of Proposition \ref{proposition3}(1):
$\Omega$ is not confined to any half-sphere.} There are several
ways of proving this.  One way is to chose $\epsilon_k$ so that
the total volume of all $\epsilon$-tubular neighborhoods is less
than half the measure of $S.$ Another is to use the fact that
$\mathcal{V}$ is the set of co-dimensional 1 subspaces normal to
edges of $P.$ Perhaps the easiest way is to choose $n$ linearly
independent vectors $v_1, \dots, v_N \in \Omega$ and observe that
$-v_1, \dots, -v_N$ are also in $\Omega.$ \qed

\bigskip
Having completed the proof of Proposition \ref{proposition3} we
have finished assembling all the elements that prove Theorem
\ref{T:main}.

%section
%section
\section{Additional Properties}
%section
%section

\begin{theorem} \label{T:byproduct1} $\empty$
 The constructed $Q$ is globally absorbing for the family
 $\{\phi_\gamma| \gamma \in P_0\}$ for any $P_0$ a compact subset of the
 interior of $P$: \textit{i.e.}, given a compact subset $X$ of $\mathbb{R}^N$
 there exists an $n_0$ (depending on $P_0$) such that if $n>n_0$ then
 $\phi_{\gamma_1}\cdots\phi_{\gamma_n}(X)\subset Q $
 for $\gamma_1,\dots,\gamma_n \in P_0.$
\end{theorem}

\medskip

We first remark that for a fixed compact subset $P_0$ of the
interior of $P$ any large enough sphere $S$ centered at the origin
will be a globally absorbing region: because any larger sphere
contracts by a fixed amount under the action of any $\phi_\gamma$,
the amount of contraction depending on distance from $P_0$ to
$\partial P.$ However, the size of the sphere $S$ varies with
$P_0.$ Theorem \ref{T:byproduct1} is a stronger result. It states
that $Q$ can be so chosen that it is universally globally
absorbing for all such $P_0.$

\medskip
\begin{proof}
Let $\rho_0> R$ where $R$ is given by \eqref{E:R}. We shall show
that the set $Q=\rho_0 Q_\infty$ is globally absorbing. From the
strengthening of \eqref{E:Q} to strict inequality which follows
from a strengthening of \eqref{E:inequality2} to strict
inequality, the other participating inequality
\eqref{E:inequality3} already being strict, we get that for any
$\gamma\in P_0$ and $\rho\ge\rho_0$ the set $\rho Q_\infty$ is
mapped by $\phi_\gamma$ into $\rho' Q_\infty$ where $\rho -
\rho'=\delta(\rho)>0 $ is a quantity which depends on the distance
from $P_0$ to $\partial P.$ The function $\delta(\rho)$ is
increasing for $\rho\ge\rho_0.$ Thus $\rho'<\rho-\delta(\rho_0)$
which leads to
\[
 \phi_{\gamma_1}\cdots\phi_{\gamma_n}(\rho Q)\subset (\rho-n\delta(\rho_0))Q
 \subset \rho_0 Q
\]
for some $n$. $\Box$
\end{proof}

\begin{theorem} \label{T:byproduct2} $\empty$
 Given a finite set of polytopes and their time dependent dynamical systems
 as defined in Theorem \ref{T:main}, there exists a $Q_\infty$
 such that the set $Q=\rho Q_\infty$ is invariant
 for large enough $\rho$ for each dynamical system.
\end{theorem}
 \medskip
\begin{proof} By virtue of the remark in \S 5.5.1 we just take $\mathcal{V}$
to contain all the co-dimension 1 subspaces orthogonal to all the
1 dimensional edges of all the polytopes.
\end{proof}

\end{document}